\documentclass[a4paper,USenglish,cleveref, autoref]{compgeomtop}

\usepackage{amsmath}
\usepackage{hyperref}
\usepackage{graphicx}
\newif\ifminted
\mintedtrue
\ifminted
\usepackage[%
frozencache=true
,cachedir=JNUMPSLA-minted-cache
]{minted}
\setminted{
  linenos=true,
  python3=true,
    autogobble,
}
\fi
\usepackage{caption}

\DOIPrefix{} %% HACK

\nolinenumbers

\makeatletter
\def\copyrightline{%
      \setbox\@tempboxa\hbox{\IfFileExists{cgt-logo.pdf}{\includegraphics[height=\EventLogoHeight\p@]{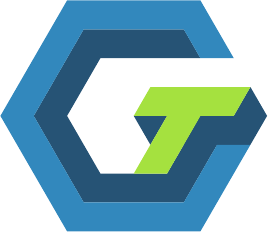}}{\includegraphics[height=\EventLogoHeight\p@]{example-image-plain}}}%
      \rlap{\hspace\textwidth\hspace{-\wd\@tempboxa}\hspace{\z@}%
            \vtop to\z@{\vskip-0mm\unhbox\@tempboxa\vss}}%
    \scriptsize
    \vtop{\hsize\textwidth
      \nobreakspace\par
      \@Copyright
      \href{https://www.cgt-journal.org/}%
      {Computing in Geometry and Topology}
      (submitted, \today) %: Volume\,\@Volume;
           pp.\,%\@ArticleNo:
           \thepage--%\@ArticleNo:
           \number\numexpr\getpagerefnumber{TotPages}\\
           }%
  }
\makeatother

\newcommand{\dts}{\mathinner{\ldotp\ldotp}}
\newcommand\psl{pseudoline}
\newcommand\psla{pseudoline arrangement}

\title{NumPSLA --- An experimental research tool for
pseudoline arrangements and order types}
\titlerunning{NumPSLA --- An experimental research tool for
  discrete geometry}

\authorrunning{G\"unter Rote}
\author{G\"unter Rote}{Freie Universit\"at Berlin, Institut f\"ur Informatik, Takustra\ss e 9, 14195 Berlin, Germany}
 {rote@inf.fu-berlin.de}{https://orcid.org/0000-0002-0351-5945}{}

\relatedversion%details%[%linktext={opt. text shown instead of the URL},
{A preliminary version of this work has been presented
  at the
  41st European Workshop on Computational Geometry
(EuroCG'25) in Liblice,
Czech Republic, April 9--11, 2025~\cite{r-nert-25}.}

\supplementdetails[%linktext={opt. text shown instead of the URL},
                   cite=github,
]{Software and computational results}{https://github.com/guenterrote/NumPSLA} %linktext, cite, and subcategory are optional

\Copyright{G\"unter Rote}

\keywords{Mathematical software, enumeration and counting, literate programming}

\acknowledgements{I thank the High-Performance-Computing Service of FUB-IT, Freie
Universität Berlin \cite{CURTA} for computing time. I thank Akiva
Weinberger (Brandeis University) for pointing out erroneous
entries in the last three lines
of
\autoref{tab:mirror} in a previous version, which uncovered a
mistake in the exclude-file for 11 points (\autoref{sec:realizable}).}
\Volume{42}
\ArticleNo{23}
\begin{document}
\maketitle
\let\paragraph=\subparagraph
  \def\subsectionautorefname{Section} % hack to avoid ssubsub... reference

\begin{abstract}
  We present a program for enumerating
all pseudoline arrangements with a 
small number of pseudolines
and abstract order types of small point sets.
This program
supports computer experiments with these structures,
and it
complements the
order-type database of Aichholzer, Aurenhammer, and Krasser.
This system makes it practical to explore the abstract order types for
12 points, and the pseudoline arrangements of 11 pseudolines.
We also report some statistics about abstract order types of 13 points
that were obtained on a compute-cluster.
\end{abstract}

\begingroup  
  \def\section*#1{}
  \bigskip
  \def\addvspace#1{\vspace {0.26pt}}

  \tableofcontents

\endgroup  
\section{Introduction}
Questions about finite configurations of points or lines are at the
core of discrete geometry.
As one example of an outstanding open question, we mention the rectilinear crossing
number problem for the complete graph $K_n$:
For a given set $S$ of $n$
points in the plane,
draw all straight segments between pairs of points in~$S$,
and count the pairs of segments that cross. What is
the smallest number that can be obtained?

\paragraph{The order type of a point set.}
The above question and 
many other questions and algorithms in discrete and computational geometry
depend only on the ``combinatorial structure'', which is typically
captured by an orientation predicate:
Consider a finite point set $S=\{p_1,\ldots,p_n\}$.
For each triplet $p_i,p_j,p_k\in S$, we need to know whether
they lie in clockwise or counterclockwise order, or whether they are collinear.
This information is enough to determine, say, the number of convex
hull vertices, or the number of  crossings.

\paragraph{The order-type database.}
It is useful if one can let the computer
exhaustively check small examples. This may provide a sanity check for
conjectures, or it may form the basis for quantitative results
that hold in general. We will mention one example below.
The prime tool that facilitates this approach
is the
order-type database of
Aichholzer, Aurenhammer, and Krasser
\cite{aak-eotsp-01,aak-eotsp-02}
at Graz University of Technology
from the early 2000's.
Originally, it contained a point set (given by coordinates)
for each of the 14,309,547 order types of 10 points (and also for the
smaller sets).
These point sets are
optimized to avoid degeneracies as much as possible.
Later, the database was extended \cite{AICHHOLZER20072} to include the 2.3 billion order types of 11
points
(see the second column of Table~\ref{tab:number}).

Over the years, the database has been enriched with all sorts of useful information about each order type,
ranging from
such basic data as the size of the convex hull to
advanced characteristics that are hard to compute, such as the
number of triangulations or the
number of crossing-free Hamilton cycles.
The  database of order types with up to 10 points can be obtained from the
website of the project
\cite{OTDB},
%\footnote
%{\url{http://www.ist.tugraz.at/aichholzer/research/rp/triangulations/ordertypes/}},
and
% the database for up to 10 points
it can be queried via an
e-mail interface.
The database for 11 points needs
102.7 GBytes (44 bytes per order type for two 16-bit coordinates per point).
Obviously, the approach of storing a representative of every order
type has currently reached its limits with 11~points.
We take an alternative approach: \emph{generating} order types
from scratch. %on the fly.

\paragraph{Big results from small sets.}
We mention just one example of a result that rests on the order-type database.
Aichh%
olzer et al.~\cite[Theorem~1]
{AICHHOLZER2020105236} proved that every set $S$ of $n$ points
in general position contains
$\Omega(n\log^{4/5}n)$ convex 5-holes, i.e., 5-tuples of points in
convex position with no points of $S$ in the interior.
Harborth~\cite{Harborth1978} showed in 1978 that every set of 10
points contains a convex 5-hole. From this, one gets an immediate
lower bound of $\lfloor n/10\rfloor$ 5-holes by partitioning $S$ into
groups of size 10 by vertical lines.
Various improvements of the constant factor of this linear bound were
obtained over the years.
The superlinear bound $\Omega(n\log^{4/5}n)$ goes beyond
what can be reached by this simple technique.
Nevertheless, at the basis of its proof,
there are some structural lemmas about sets of 11 points.
These lemmas were checked with the help of a computer by exhaustive
enumeration of order types.

\subsection{Line arrangements and pseudoline arrangements}
%\paragraph{Duality between points and lines.}
The well-known duality
\begin{equation}
  \label{eq:duality}
  \text{point }(a,b) \ \longleftrightarrow \ \text{line }y=ax-b
\end{equation}
% This duality
is a bijection between points and non-vertical lines.
It swaps the role of points and lines, and it preserves
incidences
and above-below relationships. Thus, problems about points can be
translated into problems about lines and vice versa.

\begin{figure}[htb]
  \centering
  \noindent
  \hbox to \textwidth{
\includegraphics[scale=0.8]{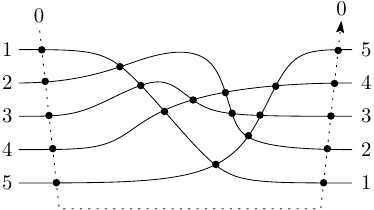}
\hfill 
\includegraphics[scale=0.85]{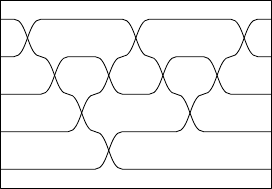}
\hfill
\hbox{\raise 0.1ex
\vbox{\small
  \catcode`\>=\active
\offinterlineskip\obeylines\openup1.2\jot\obeyspaces%
\def> #1^^M{\hbox{\texttt{#1}}}%
> 1 2-2-2-2 4-4-4-4 5
>  X       X       X
> 2 1 3-3 4 2 3-3 5 4
>    X   X   X   X
> 3-3 1 4 3-3 2 5 3-3
>      X       X
> 4-4-4 1 5-5-5 2-2-2
>        X
> 5-5-5-5 1-1-1-1-1-1
}}}
%>   1-1-1-1 4-4 5-5
%>          X   X   
%>   2 3-3 4 1 5 4-4
%>    X   X   X     
%>   3 2 4 3 5 1 3-3
%>      X   X   X   
%>   4-4 2 5 3-3 1 2
%>        X       X 
%>   5-5-5 2-2-2-2 1
  \caption{Left: A \psla\ of 5 lines, extended by a line 0 ``at
    infinity''.
  Middle and right: wiring diagrams}
  \label{fig:psa5}
\end{figure}

\paragraph{Pseudoline arrangements and abstract order types of
  points.}
Pseudoline arrangements are a generalization of line arrangements.
A \psla\ (PSLA) is a collection of unbounded curves, with the condition
that any two curves intersect exactly once, and they cross at this
intersection point.
We refer to these curves as \emph{\psl s} or simply as \emph{lines}.
See Figure~\ref{fig:psa5} for an example with 5 pseudolines.
The middle and the right picture show a standard representation
%of a PSLA
as a \emph{wiring diagram}, in two different styles. Both drawings are produced by
our program.
In a wiring diagram, the pseudolines 
run on $n$ horizontal tracks, and they cross by swapping between
adjacent tracks.

By duality, there is an analogous notion for point configurations, an
\emph{abstract order type} (AOT).
We will elucidate this relation in Section~\ref{sec:dual}. There is a variety
of equivalent notions for these objects, such as rhombus tilings, oriented matroids of rank 3, or signotopes; see for
example \cite[Chapter 6]{Felsner}.

Our program focuses on pseudoline arrangements as the primary objects.
The main reason is that they are easy to generate in an incremental way.
Another % second%ary
advantage %reason
is that they are easy to draw and to visualize.

Throughout this paper, we will assume general position. In other
words, we restrict our attention to \emph{simple} \psla s, where no
three lines go through a common point. In the setting of point sets,
this corresponds to forbidding collinear point triples.

\begin{table}[htb]
  \centering
  \vbox{
    \def\gobble#1{}
    \tabskip=0pt plus 1em
\halign to \textwidth{\strut\hfil$#$ &\ \hfil$#$ &\ \hfil$#$ &\ \hfil$#$ &\quad$#$\hfil
  &\gobble{#}%\ \hfil$#$
  &\hskip-0em\hfil$#$\cr
&\hbox{[\href{https://oeis.org/A006247}{A006247}]}
&\hbox{[\href{https://oeis.org/A063666}{A063666}]}
&\Delta = {}&&&\hbox{[\href{https://oeis.org/A006245}{A006245}]}
 \cr 
 n & \hbox{\#AOT} & \hbox{\#OT}
 & \hbox{\#nonr.\ %ealizable
   AOT}
 %\Delta
 &\smash{\raise 0.7em\hbox{$\displaystyle\frac{\Delta}{\hbox{\#AOT}}$}}
 & \hbox{ \ \#mirror-symmetric AOT}
 & \hbox{\#%$x$-monotone
   PSLA}
\cr
\noalign{\hrule}
3 & 1 & 1 & 0 & 0 & 1 & 2 \cr
4 & 2 & 2 & 0 & 0 & 2 & 8 \cr
5 & 3 & 3 & 0 & 0 & 3 & 62 \cr
6 & 16 & 16 & 0 & 0 & 12 & 908 \cr
7 & 135 & 135 & 0 & 0 & 28 & 24{,}698 \cr
8 & 3{,}315 & 3{,}315 & 0 & 0 & 225 & 1{,}232{,}944 \cr
9 & 158{,}830 & 158{,}817 & 13 & 0.01\,\% & 825 & 112{,}018{,}190 \cr
10 & 14{,}320{,}182 & 14{,}309{,}547 & 10{,}635 & 0.07\,\% & 13{,}103 & 18{,}410{,}581{,}880 \cr
11 & 2{,}343{,}203{,}071 & 2{,}334{,}512{,}907 & 8{,}690{,}164 & 0.37\,\% & 76{,}188 & 5{,}449{,}192{,}389{,}984 \cr
12 & 691{,}470{,}685{,}682 &  &  &  &  &
2{,}894{,}710{,}651{,}370{,}536 \cr
13 & 366{,}477{,}801{,}792{,}538
&  &  &  &  &
2{,}752{,}596{,}959{,}306{,}389{,}652\cr
%4675651520558571537540
}}

\smallskip
\caption{\#AOT = number of abstract order types for $n$ points.
  \#OT = number of order types.
  \#PSLA = number of ($x$-monotone) pseudoline arrangements with $n$ pseudolines.
  %$\Delta$ refers to the number or
  % fraction of nonrealizable AOTs.
These are the objects that the program actually
enumerates
one by one  (almost, because we try to apply shortcuts).
  The column headings link to the
    corresponding entries of the
Online Encyclopedia of Integer
Sequences \cite{OEIS}.
%These numbers are known up to $n=15$.
  }
  \label{tab:number}
\end{table}

\subsection{Overview}
We will describe our algorithm for enumerating \psla s, and we will
apply it to enumerate abstract order types.
None of the techniques that we use are novel, but we have tried to
streamline and simplify the algorithms.
In terms of speed, we can compete with the order type database
(see
\autoref{sec:benchmark} for a direct comparison).
The main distinction is, of course, that the order type database contains only
\emph{realizable} order types, and that they come with coordinates.
For many applications, the restriction to
{realizable} order types is not important, and coordinates are
not needed. In those applications, our approach shows its strength.
Mustering the 14 million 10-point abstract order-types takes 10--20
seconds.
The 11-point sets can be handled in half an hour, % of CPU time.
and the 12-point sets take about 200 CPU hours.
To this, one must of course add the time for whatever one
wants to do with those order types.
The program is trivially parallelizable, and with a powerful
compute-cluster, it is feasible to go even for 13 points,
see Section~\ref{some-results}.

The program is available in a GitHub repository~\cite{github}.
%at \url{https://github.com/guenterrote/NumPSLA}.
It is written in the programming language~C, using the CWEB system of structured
documentation of Donald E. Knuth and Silvio Levy\footnote{\url
  {http://tug.ctan.org/info/knuth/cwebman.pdf}}.

We have occasionally used the enumeration
for our research questions
\cite{bckr-sxmp-25},
and we hope it finds other users.
%and it been used by
%other researchers as well
%\cite{parpalak2026, savchuk2025}.

%We encourage everybody to try out whether they can
%use it for their own purposes.

\section{Enumeration of \psla s}
\label{sec:Enumeration}
We concentrate on $x$-monotone \psla s, in which the curves are
$x$-monotone.
Every \psla\ can be drawn in an $x$-monotone way, but this incurs a
choice: One of the unbounded faces must be selected as the
\emph{top face} $T$, and the opposite unbounded face will become the
\emph{bottom face} $B$.
Then the lines run from left to right, and we number them from 1 to
$n$
as they appear from top to bottom on the left side (see Figure~\ref{fig:psa5}).
If they were straight lines, they would be numbered by increasing
slope.

\subsection{Representing a \psla}
\label{sec:pred-succ}
\begin{figure}[htb]
  \centering
  \includegraphics{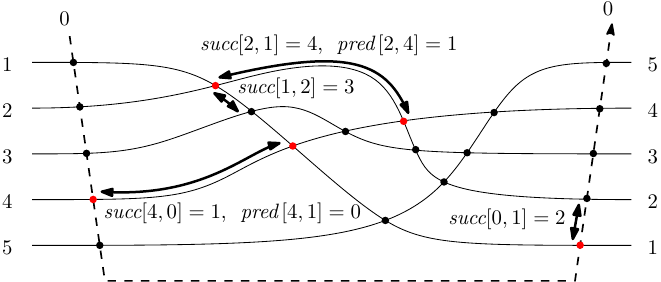}
  \caption{Some successor and predecessor pointers for the PSLA of Figure~\ref{fig:psa5}}
  \label{fig:successors}
\end{figure}

The vertices and edges of a \psla\ form a planar graph. The storage
and manipulation of this graph is greatly simplified by the fact that
we have a precise control over the vertices: There is a vertex for
each pair of lines, and every vertex has degree~4.
We thus store the edges in two 2-dimensional arrays
\textit{succ} 
and \textit{pred} of successor and predecessor pointers,
see \autoref{fig:successors} for an illustration.
The entries
$\textit{succ}[j,k]$
and $\textit{pred}[j,k]$
refer to the
crossing between line $k$ and the line~$j$.
We think of the lines as oriented from left to right.
Then
$\textit{succ}[j,k]$ and $\textit{pred}[j,k]$ point to the next and
previous crossing on line~$j$.
For the reversed index pair $[k,j]$, we get the corresponding information for line $k$.
Thus, in
 \autoref{fig:successors},
%the example of Figure~\ref{fig:psa5},
%$\textit{succ}[2,3]=5$,
% and accordingly, $\textit{pred}[2,5]=3$.
$\textit{succ}[2,1]=4$,  and accordingly,
$\textit{pred}\,[2,4]=1$.

 We can easily determine which of $j$ and $k$ enters the intersection
 $(k,j)$ from the top and bottom: By our numbering convention, the
 line with the smaller index always enters above the other line, and
 to the right of the crossing, it lies below the other line.
 
The infinite rays on line $j$ are represented by the additional line 0:
$\textit{succ}[j,0]$ is the first (leftmost) crossing on line $j$,
and
$\textit{pred}[j,0]$ is the last crossing.
The intersections on line 0 are cyclically ordered $1,\dts,n$.
Thus, $\textit{succ}[0,i]=i+1$ and $\textit{succ}[0,n]=1$.

\begin{figure}[htb]
  \centering
  
  \includegraphics{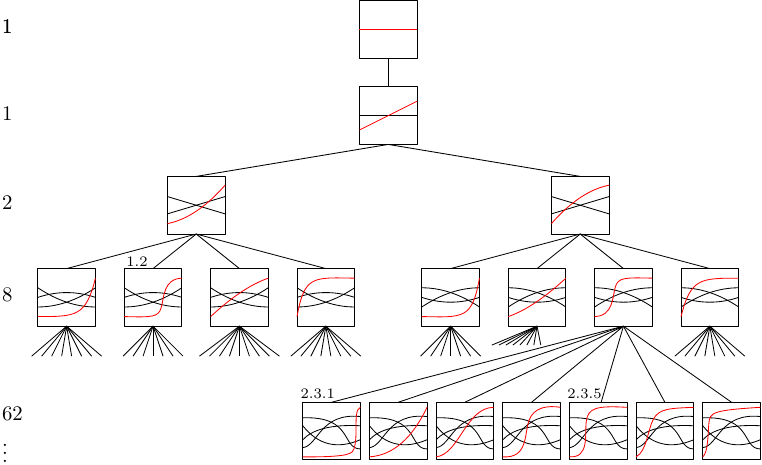}
  
  \caption{The first three levels of the enumeration tree
    and a few nodes of the fourth level. The last inserted \psl\ is
    highlighted in color. For selected nodes, the Dewey decimal notation
    is given.
%    and accession number-WEG!
  }
  \label{fig:tree}
\end{figure}

\subsection{Incremental generation of \psla s}
\label{sec:incremental}
We generate a PSLA with $n$ lines by inserting line $n$ into a
PSLA with $n-1$ lines, in all possible ways.
Then each PSLA has a unique predecessor PSLA, and this
predecessor relation
imposes a tree
structure on the PSLAs, as shown in Figure~\ref{fig:tree}. Our program
explores this tree in depth-first order.
If we number the children of each node in the order in which they
are visited, this leads to
a
unique
identifier for every node, and thus for every PSLA,
analogous to the Dewey decimal classification that is used to classify
books in libraries.

\begin{figure}[htb]
  \centering
  \includegraphics{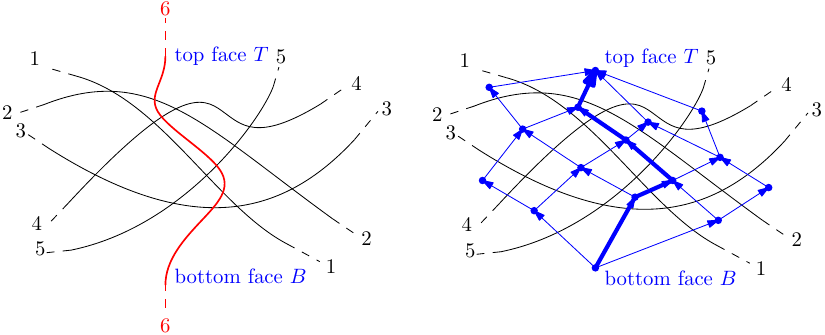}
  \caption{Left: Threading line 6 through a PSLA of 5 lines.
    Right: The dual DAG of this PSLA}
  \label{fig:dualDAG}
\end{figure}

Inserting the $n$-th pseudoline into a PSLA of $n-1$ lines corresponds to
threading a curve from the bottom face $B$ to the top face $T$,
see Figure~\ref{fig:dualDAG}.
(We temporarily relax the requirement that the extra \psl\ has to be
$x$-monotone.) Following Knuth~\cite[Section~9, p.~38]{axioms-hulls},
such a curve is called a
\emph{cutpath}
\cite{felsner-valtr-11,d-ibnpa-25}. This corresponds to a source-to-target path in the dual
graph of the PSLA. Orienting the dual edges in the way how line $n$
can cross them, namely, from below to above, leads to a directed
acyclic graph (a DAG).
We can enumerate all such paths in a backtracking manner.
Since the DAG has no sinks other than the target vertex $T$, a path
cannot get stuck, and thus the enumeration of the paths is simple and fast.

\begin{figure}[htb]
  \centering
  \includegraphics{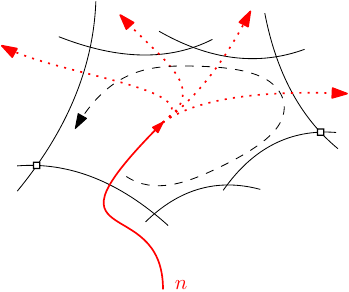}
  \caption{Continuing line $n$ after entering a face.}
  \label{fig:explore-face}
\end{figure}

The whole algorithm is thus a double recursion.
The outer recursion extends a PSLA by adding a \psl\ $n$.
The inner
recursion extends a partially drawn pseudoline $n$ to the next crossing,
see Figure~\ref{fig:explore-face}.
It is implemented by walking along the boundary of the face that has
been entered through the last crossing.
All upper edges of the face are candidate edges for the next crossing
of line $n$, and we try them in succession.
We have decided to walk in counterclockwise order around the
face. This means that the paths for line $n$ are generated in
``lexicographic'' order from right to left,
as can be checked in Figure~\ref{fig:tree}.

Appendix~\ref{python} shows a self-contained \textsc{Python} program
that implements this enumeration algorithm.

\section{Duality between \psla s and abstract order types}
\label{sec:dual}
\begin{figure}[htb]
  \centering
  \hbox to \linewidth{\hss
  \includegraphics[scale=0.9]{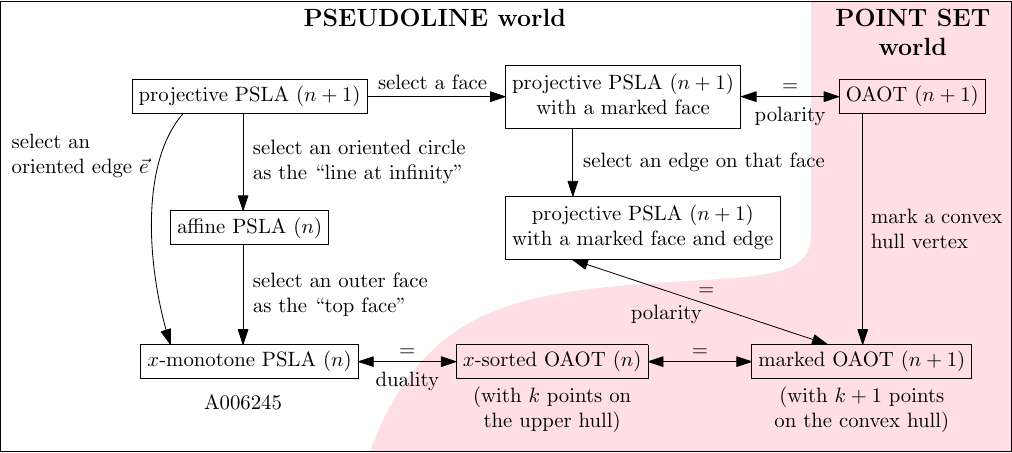}\hss}
\caption{Relation between different concepts.
 An arrow in one
  direction
  indicates a specialization.}
  \label{fig:PSLA}
\end{figure}

The
duality between \psla s and abstract order types is not as
straightforward as one would hope for. Figure~\ref{fig:PSLA} 
shows the intricate network of relationships.
The top right box refers to \emph{oriented} abstract
 order types (OAOTs), where a point set is still distinguished from its
 reflection. (AOTs don't make this distinction.)
At the lower left corner, we find our favorite objects, the
($x$-monotone) PSLAs. %with $n$ pseudolines.

There is an offset of 1 in the number of objects when transferring
between the ``point set world'' (shaded) and the ``pseudoline world'':
The
($x$-monotone) PSLAs
with $n$ pseudolines correspond to OAOTs or AOTs with $n+1$ points.
However, the correspondence is not one-to-one.
Different  PSLAs
may give rise to the same OAOT and AOT, and the algorithm has to take
care of this ambiguity when trying to enumerate OAOTs or AOTs without
duplication. The details are given in 
\autoref{sec:duplicates}.

\begin{figure}[htb]
  \centering
  \includegraphics{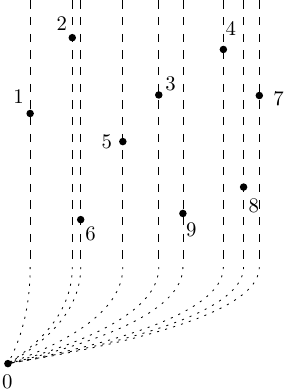}
  \caption{Modeling the
    $x$-sorted order of a point set by an extra
    point $0$.}
  \label{fig:OT-points}
\end{figure}

%In the lower left corner of
%Figure~\ref{fig:PSLA}
%(p.~\pageref{fig:PSLA}), we find our favorite objects, the
Let us discuss some of the boxes in more detail,
starting at 
 the lower left corner of~Figure~\ref{fig:PSLA}:
The \psl s of an $x$-monotone PSLA
are numbered from 1 to $n$ in order of increasing slope.
If we start with the analogy of a line
arrangement and apply
the duality~\eqref{eq:duality}, we get a set of points that
are sorted by $x$-coordinate,
as shown in Figure~\ref{fig:OT-points}.
Now, the notion of being sorted by $x$-coordinate is foreign to order
types, but we can incorporate it by imagining a point $0$ at vertical
positive infinity, around which the points are radially sorted.

We can also move this extra point to a finite distance, sufficiently
far above, without changing the order type. Moreover, if we move the
point $0$ to
the left of all points, as
 in Figure~\ref{fig:OT-points}, we see that we can let this point
 correspond to the line 0 in the PSLA, or more precisely, to the part
 of line 0 that lies at the left of all crossings.
 This line has a smaller slope than all other lines, and
 it intersects the other lines in the order $1,\ldots,n$.
 This is perfectly in accordance with the situation in the dual setting:
 The corresponding point $0$ has
 a smaller $x$-coordinate than all other points, and the counterclockwise cyclic
 order of the other points is $1,\ldots,n$.

 This extended point set has $n+1$ points, and it has a special point $0$
 on the boundary.
 We see that this is equivalent to an arbitrary set of $n+1$ points where some
\emph{pivot point} on the convex hull is marked. By a projective transformation,
 the pivot point can be moved far up without changing the order
 type. %, the set can be rotated to make the direction almost vertical.

 Thus we have explained the three boxes in the bottom row of
 Figure~\ref{fig:PSLA}. The boxes refer to \emph{oriented} abstract
 order types (OAOTs), because at this point, we still distinguish a point set from its
 reflection.

\begin{figure}[htb]
  \centering
  \includegraphics{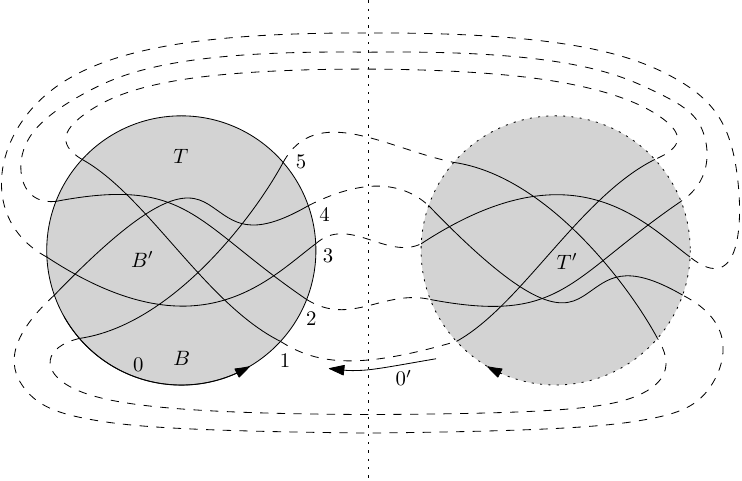}
  \caption{The spherical model of a projective PSLA,
for the PSLA of Figure~\ref{fig:dualDAG}
  }
  \label{fig:spherical-model}
\end{figure}

 We are, however, interested in point sets without a marked pivot
 point.  Therefore we must study what it means to select another hull
 point as the pivot point.  This is best understood by looking at
 \psla s in the projective plane.  As the model of the projective
 plane, we use the sphere in which opposite points are
 identified. Figure~\ref{fig:spherical-model} shows a picture, where
 the image of the PSLA to which we are used appears on the ``front half''
 of the sphere in the left part of the picture. The ``back half'' of
 the sphere, which carries the centrally reflected PSLA, is unfolded
 into the right part of the picture, so that we look at both parts
 from the outside.  On the sphere, each \psl\ becomes a closed cycle.
 The line 0 ``at infinity'' is the cycle that separates the front part
 from the back part.  The dashed lines connect points that are
 identical, as the front part and the back part are stitched together.
 %(The figure is
 %geometrically inaccurate insofar as each curve ought to cross the
 %circle 0 at antipodal points. The front surface and the back surface
 %would have to be deformed when being stitched together.)

Now, on this spherical model,
 we have $n+1$ lines. They are
all equal; line 0 does not
play a distinguished role.
In fact, the
\textit{succ} 
and \textit{pred} pointers allow navigation on the sphere
just fine. If we follow
the \textit{succ} pointers along some line $j$ without caring to stop
when we cross line $0$, we will simply traverse the whole cycle again
and again.

We have obtained our $x$-monotone \psla\ because we know which circle
is line 0, and moreover, we have marked two opposite faces within this
circle as the bottom face $B$ and the top face $T$.

We can obtain another $x$-monotone \psla\ from the same projective
class
by declaring a different line to be line 0, and marking the faces that
should become the bottom and top faces.
One such choice is indicated by the labels $0',B',T'$ in 
Figure~\ref{fig:spherical-model}.
A different way to express this is to say that we pick a directed
edge as a \emph{starting edge}, namely the edge of cycle $0$ that has
the face $T$ on its left.

This discussion covers the boxes on the left side of
Figure~\ref{fig:PSLA}.
As an intermediate notion, we have \emph{affine} (or \emph{Euclidean})
PSLAs, where the line at infinity is fixed, but it has not been
decided which unbounded faces are the bottom and the top faces.
%The boxes in Figure~\ref{fig:PSLA} refer to \emph{oriented} abstract
%order types (OAOTs), because at this point, we still distinguish a point set from its
%reflection.
Figure~\ref{fig:PSLA} includes some intermediate boxes, in which some
data are partially fixed, and it shows a few more translations between the \psl\
world and the point-set world, but we don't discuss them here.

If a different starting edge has been chosen, it is not hard to
realize this in the data structure. The graph is the same as before;
one just has
to relabel the lines. The line through the starting edge becomes line $0$, and
the other lines get the labels $1,2,\ldots,n$ in the order in which
they are crossed by line 0, starting from the starting edge.
We simply need to carry out this relabeling for $j$, $k$, and $i$ or $i'$
in all relations
$\textit{succ}[j,k]=i$
and $\textit{pred}[j,k]=i'$.

\subsection{The convex hull in the \psl\ world}
\begin{figure}[htb]
  \centering
  \includegraphics{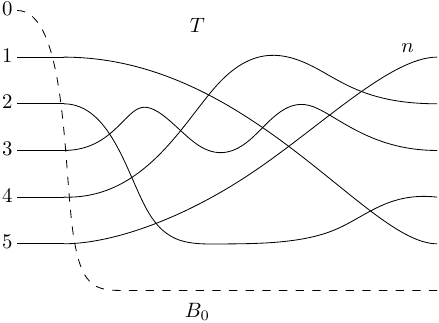}
  \caption{The convex hull in a \psla}
  \label{fig:PSLA-hull}
\end{figure}

As is well-known, the convex hull of a point set consists of those
points whose dual line is incident to the top face or the bottom face.
However, when applying this criterion, we must add line 0 as the line
with the most negative slope, as illustrated in
Figure~\ref{fig:PSLA-hull}.
Then there are only two lines incident to the bottom face: lines 0 and
$n$.
But these two lines are anyway also incident to the top face. Thus,
in our setting, the convex hull vertices correspond to the edges of
the top face (including edge 0).

\subsection{The orientation predicate}
  %.Working with an abstract order type}
\label{sec:working}

\begin{figure}[htb]
  \centering
  \includegraphics{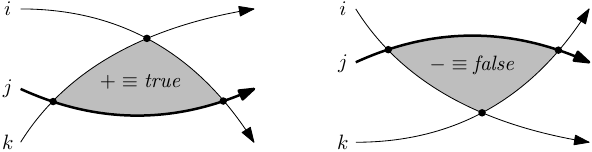}
  \caption{The orientation of three lines}
  \label{fig:orientation-pred}
\end{figure}

The 
\textit{succ}
and \textit{pred} arrays are useful for navigating in the arrangement,
but to get the full power of working
with an abstract order type, one needs the orientation predicate.
In terms of \psl s, the orientation is defined as shown in
Figure~\ref{fig:orientation-pred}. For three lines $i<j<k$ the
orientation is determined by looking at the triangle formed
by these lines. The orientation
\textit{orient}$(i,j,k)$
is positive if the triangle lies below
$j$ and negative otherwise.
The orientation is unchanged under an even permutation of the parameters
$(i,j,k)$, and it is flipped by an odd permutation of the parameters
$(i,j,k)$.
This orientation agrees with the orientation of the
corresponding point set, in case we apply the duality~\eqref{eq:duality} to a proper line
arrangement, in the following sense:
For three lines $i<j<k$ with dual points $p_i,p_j,p_k$,
the crossing of $i$ and $k$ lies below $j$ if and only if
$p_j$ lies below the segment $p_ip_k$, i.e., if and only if
$p_i,p_j,p_k$ are in counterclockwise order.
(For triples involving line~$0$, this requires that the dual point~$0$
lies far above and to the left of the other points, see below.)
Extending the definition to \psla s is in fact one way to define
abstract order types.

Now, for $i<j<k$, as shown in the picture, the orientation can be
figured out if one knows the order of the crossings
along line $j$, for example: Is the crossing $(j,i)$ to the left or to
the right of $(j,k)$?
This information is not available, but it can easily be provided by
preprocessing.

Thus, when we want to work with a PSLA, we prepare additional data
structures,
the \emph{local sequences array} $P$ and the
\emph{inverse local sequences array} $\bar P$.

\subsection{The local sequences matrix and its inverse}
Here is a representation as a two-dimensional array.
For
each pseudoline $i$, the sequence $P_i$ indicates the sequence of crossings with the other
lines, starting at $0$ % by convention
and moving to the right.
For the example in Figure~\ref{fig:psa5}, the local sequences are as follows:
\begin{displaymath}
\vbox{
\def\minus{{\setbox0=\hbox{0}\hbox to \wd0{\hss-\hss}}}
\catcode`-=\active\let-=\minus
\begin{tabbing}
  \qquad\=\+
  $P_0=[1,2,3,4,5]$\qquad\qquad\=$\bar P_0=[-,0,1,2,3,4]$
  \qquad\qquad\=%$B_0=[0,0,0,0,0]$ 
  \\$P_1=[0,2,3,4,5]$\> $\bar P_1=[0,-,1,2,3,4]$ %\>$B_1=[0,0,0,0,0]$ 
  \\$P_2=[0,1,4,3,5]$\> $\bar P_2=[0,1,-,3,2,4]$ %\>$B_2=[0,0,0,0,1]$  
  \\$P_3=[0,1,4,2,5]$\> $\bar P_3=[0,1,3,-,2,4]$ %\>$B_3=[0,1,0,0,1]$ 
  \\$P_4=[0,1,3,2,5]$\> $\bar P_4=[0,1,3,2,-,4]$ %\>$B_4=[0,1,1,1,0]$ 
  \\$P_5=[0,1,2,3,4]$\> $\bar P_5=[0,1,2,3,4,-]$ %\>$B_5=[0,1,1,1,1]$  
\end{tabbing}
}
\end{displaymath}
The first row $P_0$ and the first column are fixed. % determined.
Each row $P_i$ consists of $n$ different elements, excluding the
element $i$ itself.
The inverse local sequence
$\bar P_i$ is essentially the inverse permutation of $P_i$:
The $j$-th element of $\bar P_i$ gives the position in $P_i$ where the
crossing with $j$ occurs, where positions are counted from~$0$.
The diagonal entries are irrelevant.
From the \textit{succ} links, it is straightforward to build
the local
sequences
and the inverse local
sequences,
in $O(n^2)$ time.
With the help of $\bar P$, the orientation predicate can be evaluated
in constant time as the exclusive-or of three simple tests:
%and the reverse local
%sequences,
\begin{displaymath}
\textit{orient}(i,j,k) \equiv (i<j) \oplus (j<k) \oplus (\bar P[j,i]<\bar P[j,k])
\end{displaymath}
Here we use a Boolean value instead of a sign~$\pm$.
It is clear that this formula is correct for the standard case
$i<j<k$, and
it is easy (but tedious) to check
% we have confirmed
that it works for all other orderings of $i,j,k$.

\section{Elimination of duplicate AOTs}
\label{sec:duplicates}

An abstract order type with $h$ hull vertices corresponds to $2h$
PSLAs:  Each of the $h$ hull vertices can become vertex~0, and there are two choices of
orientation.
(If there are symmetries, some of these $2h$ PSLAs will coincide.)
The standard approach to tackle this problem is to compute some sort of
canonical representation.
Since a \psla\ is a plane graph, one may apply isomorphism testing
algorithms for plane graphs.

Indeed, such methods are available, due to
Hopcroft and Wong
\cite{hopcroft+wong-isomorphism-74} from 1974,
for (3-connected) planar graphs, and
a simplified method of
Akutsu, de la Higuera, and Tamura~\cite{akutsu_CPM.2018}
for plane graphs with a fixed outer face.
Both methods run in linear time in the size of the graph,
which is $O(n^2)$ in our case, and they
 produce a canonical
 labeling as well as the automorphism group of the graph.
We have to adapt the input, because \psla s have infinite rays.
A simple way to treat these rays is to enclose the arrangement by a
ring and clip the rays at the ring, generating a cycle of $2n$ degree-3 vertices.

(We mention another special algorithm that tests
(abstract) \emph{order-type} isomorphism of \emph{point sets}
and determines the automorphism group in $O(n^2)$ time~\cite{order-type-isomorphism}; this algorithm attacks the problem directly and
does not rely on %have recourse to
planar graph isomorphism via duality.)

However, instead of applying these general techniques for planar
graphs,
we apply %have chosen
a simpler, tailored method, which
runs very fast in practice.
To eliminate duplicates, we compare the
local sequences matrices
$P$ lexicographically. The algorithm produces an AOT $A$ only if
the $P$-matrix of the PSLA at hand is the smallest in the
class of PSLAs that represent~$A$.
Conceptually, we proceed as follows:
Let $A^1=A,A^2,A^3,\ldots,A^{2h}$ be the $2h$ ways of labeling
the point set $A$, preserving the set of convex hull vertices.
We then look at the corresponding local sequences matrices
$P^1,P^2,\ldots,P^{2h}$.
 $P^1$ represents the \psla\ currently at hand, and
%look at the current local sequences matrix
$P^2,\ldots,P^{2h}$ are its competitors.
If the current matrix is not the smallest, we discard the current
PSLA.
On this occasion, we will also find out when some of the other
$P$-matrices are equal to $P^1$. This indicates the presence of a symmetry.
The symmetry may be a rotational symmetry, rotating the convex $h$-gon
by some number of vertices.
(The first $P^j$ that is equal to $P^1$ determines the period, which must be a divisor of~$h$.)
 A~mirror symmetry can occur alone or in combination with a
 rotational symmetry, and it will also be detected.

 There are several considerations that play a role
 in practice:
 \begin{enumerate}
 \item 
 The average number $h$ of sides of the convex hull is a
 little bit less than 4, see~Table~\ref{tab:hull}. This confirms
 theoretical predictions of 
 Goaoc and Welzl
  \cite{goaoc-welzl-2023}.
 They showed
 that for
\emph{labeled} order types (where symmetries don't matter), the
average
size of the convex hull
is
\begin{equation}
  \label{eq:welz}
4-\frac8{n^2-n+2}.
\end{equation}
This statement holds both for
AOTs
 \cite[Theorem 10.2]{goaoc-welzl-2023}
and for
(realizable) \emph{order types}
\cite[Theorem 1.2]{goaoc-welzl-2023}.
In the latter setting of order types, convergence to 4 carries over to the unlabeled case
\cite[Theorem 1.3]{goaoc-welzl-2023}.
In our setting of unlabeled abstract order types, no such convergence
result has been proved. Nevertheless,
Table~\ref{tab:hull} shows that formula
    \eqref{eq:welz} seems to give a very precise estimate even in this setting.

  \item The vast majority of AOTs have no symmetries.
    (See
\autoref{tab:symmetries} in
    \autoref{sec:symmetries} for exact counts of symmetric versus
asymmetric  AOTs.)
   Thus we can assume that roughly only one out of 8 PSLAs
  is the lex-min PSLA, and 7 out of 8 are generated in vain.
  One can check that the numbers of Table~\ref{tab:number} are in
  accordance with this.
  The 112,018,190 PSLAs with 9 lines give rise to only
  14,320,182 AOTs with 10 points.
  The ratio is 0.127838, just barely larger than 1/8.
  
\item Most of the total runtime for enumeration is spent in the lex-min test at the leaves
  of the tree.
 \end{enumerate}

 In practice, it would be wasteful to compute the complete matrices
 $P^1,\ldots,P^{2h}$ in advance, which would take $\Theta(hn^2)$ time. We compute the first entry of each
 matrix and compare these entries. It may turn out at this point that $P^1$ has
 already lost, and we can quickly abandon the comparison. Some other
 matrices might also be out of the game, and they are discarded.
 For the matrices that remain, we compute the second entry, and so on
 (see also
 \cite[p.~4]{AICHHOLZER20072}).
 The comparison will only go to the very end if some matrices are
 equal,
 and this can only happen in case of symmetry. As mentioned, symmetric
 solutions are a tiny fraction. %small minority.

 \subsection{Early screening of \psla s}
 The way we compare the local sequences matrices in the lexicographic order is row-wise from
 right to left. That is, we start with the right-most entry $P_{1n}$
 in the first row $P_1$. (Row $P_0$ is always the same.
 In contrast to the convention for $\bar P$, we now count the columns~$j$
 of the entries $P_{1j}$ from 1 to~$n$.)
 The reason for this unusual choice is that, in some preliminary tests, it seemed to be more effective in
 connection with the screening approach that is described below.

 We have mentioned that the effect of choosing a different starting edge
 consists of relabeling all lines. Thus,
 in order to compute the matrices
 $P^2,\ldots,P^{2h}$, we compute a renaming table for
 each
 matrix. This takes $O(n)$ time per matrix, by simply following the
 pointers along one pseudoline.
  This task has to be completed before the first matrix entry is even
  looked at.

To speed things up, we sidestep the renaming table and compute the entry
 $P_{1n}$ (and only this entry) directly. The meaning of $P_{1n}$ is the (label of) the last
 line $\ell$
intersected by line~1. This label is defined by how far away the
intersection of $\ell$ % this line
is from the start, when walking along
line~0.
This interpretation can be used to determine the value of $P_{1n}$
even with incorrect labels,
by simply walking along
line~0 (in a \emph{pedestrian} way, so-to-speak).

If, for example, one of the other matrices $P^i$ has a
smaller value
$P_{1n}^i$ than
$P_{1n}^1$
in the matrix $P^1$, we immediately conclude that $P^1$ is not lex-min, and we have saved a lot of work.
A matrix $P^i$ with
$P_{1n}^i>P_{1n}^1$ can be excluded from further consideration,
and hence its renaming table need not be computed.
 The details are a bit tricky, and they are explained in the
 documentation of the program.
 This screening test is quite effective.
 For example,
 there are
 18,410,581,880 PSLAs 
 with $n=10$ lines.
 Of these, only
 5,910,452,118 pass the screening test.
 Eventually, only
 2,343,203,071 PSLAs are really lex-min, and this is the number of
 AOTs that we really want.

For those cases that pass the screening test, 
it turns out that the lex-min testing procedure is quite fast:
When enumerating
 AOTs with $n\ge10$ points, on average, a lex-min test had to look at
 less than 6 entries in total before it could make a decision.
This total is over all matrices
$P^1,\ldots,P^{2h}$ taken together.
(This does not include
the $2h$ entries $P_{1n}^i$ that were compared in
the screening test. The screening test eliminates some of
the $2h$ candidates, but for the surviving candidates,
the lex-min test looks at the entries $P_{1n}^i$ again, for uniformity.
By adapting the code, the 6 entries that are looked at on average
could be further reduced.)

%, excluding the ones
%that were eliminated by

%(These data were measured in combination with
%the pre-screening procedure at level $n-1$ that is described in the next section.
%That procedure already eliminates many ``easy'' PSLAs which can be classified
%as being not lex-min by looking at the entry $P_{1n}$.)
% THAT SHOULD BE IRRELEVANT BECAUSE THOSE WOULD
% BE THROWN OUT BY THE SCREENING TEST.

% 5.056 for n=12
% 5.209 for n=11
 
% 5.699 for n=9
% 6.121 for n=8
% 7.723 for n=7

 \subsubsection{More aggressive pre-screening at the next-to-last level}
 In some cases, it can be determined already at level $n-1$ that there
 is no way that the insertion of line $n$ into the current PSLA can
 lead to a lex-min PSLA. In this case, we can abandon the
 procedure right away, instead of generating all
 children in the tree and subjecting them to the lex-min test.
The details are described in the documentation of the program.

\section{Parallelization}
 We have implemented a trivial way to parallelize the enumeration.
The user can choose a \emph{split level}, usually~8. The program will then work
 normally up to level 8 of the tree, that is, it will
 enumerate all 1,232,944 PSLAs with 8 lines, but it will only expand
 a selection of these PSLAs.
 The selection is determined as follows. As the PSLAs with 8 lines are
 enumerated, a running counter is incremented, thus %implicitly
 assigning a number
 between 1 and 1,232,944 to each PSLA. We specify a modulus $m$ and a
 value $k$. Then the program 
  will expand only those nodes whose number is congruent to $k$
  modulo~$m$.
  By running
 the program for $k=1,\ldots,m$, the work is split into $m$ roughly equal parts.
 The same parallelization technique is used in the
\texttt{gtools} 
graph generation programs that
% \texttt{geng}
% program,
% which
are part of the \texttt{nauty} system\footnote{
\url{https://users.cecs.anu.edu.au/~bdm/nauty/} or
  \url{https://pallini.di.uniroma1.it/}}, % for generating graphs,
see \cite[Section~8] {McKay1998}.

\section{Enumerating only the realizable AOTs}
\label{sec:realizable}

We have implemented a provision to enumerate only the (realizable) order
types of points sets, for up to 11 points,
to make the results comparable with those of the order-type database:
There is an option to specify an \emph{exclude-file} for the program. The
exclude-file is a sorted list of decimal codes for tree nodes that should be
skipped.\footnote{Currently the exclude-file feature does not work together with the
  parallelization feature. (For 11 points, the program should anyway be
  fast enough without parallelization.)}
The exclude-files were prepared with the help of the order-type
database.
Essentially, we are storing the AOTs that are \emph{not} realizable,
which is a tiny minority compared to the realizable ones, see
 Table~\ref{tab:number}.
Still, the exclude-file for up to 11 points has
7,710,079 entries
and needs
 161.9 MBytes.
%exclude11.txt has
%8{.}699{.}559 entries and 184{.}634{.}808 bytes
 (With some technical effort, like eliminating common prefixes
 or a compressed binary format, one could reduce this
 space requirement significantly.)
 
% one could reduce this
% To save space, one could use a compressed binary format with
%42 bits per entry, which would amount to 45.7 MBytes.)
% would save space.)
%42 bits = 1+2+4+5+6+7+8+9 bits per entry
%>>> 8699559 * 42 / 8
%45672684.75 ~ 45.7 MBytes
% text-based format with prefix compression:
% 36933905 bytes, see shared-prefixes-exclude-file-compress.py

 \section{Symmetries of \texorpdfstring{$x$}{x}-monotone PSLAs}
 \label{sec:symmetry-PSLA}

For \emph{$x$-monotone} PSLAs, there are three natural symmetry operations.
%that one could consider:
\begin{enumerate}[a)]
\item A horizontal %mirror
  reflection (with a
{vertical} mirror axis) keeps the bottom and the top
faces fixed and exchanges the
 \textit{pred}
 and \textit{succ} pointers.
 The lines
 are traversed in the
   opposite direction, and they
 are renumbered from $1,2,\ldots,n$ to $n,n-1,\ldots,1$.
(Such a vertical symmetry corresponds to a mirror symmetry of the
associated AOT, with the mirror going through the pivot point.)
\item A vertical %mirror
  reflection %-symmetry
(with a
{horizontal} symmetry axis) swaps the bottom face with the top
face.
As above, the lines are renumbered from $1,2,\ldots,n$ to $n,n-1,\ldots,1$.
\item
   A $180^\circ$ rotation is obtained as a combination of the above
   operations.
   The lines are not renumbered, but they are traversed in the
   opposite direction.
\end{enumerate}

These symmetries
are the natural symmetries for $x$-monotone PSLAs.
They preserve many properties of PSLAs, for example, the number of
cutpaths.
By taking advantage of these symmetries, one can typically reduce
the amount of computation for experiments by a factor of up to four.

A PSLA with an odd number $n\ge3$ of lines cannot have
a horizontal symmetry axis:
The middle line $(n+1)/2$ must pass either above or below the
crossing of $1$ and $n$, and this property is inverted by a reflection at a
horizontal axis.

A PSLA with $n\ge3$ lines cannot have a $180^\circ$ rotation.
For $n=3$ (and general odd $n$), the same argument as above
can be applied.
%\emph{both} a vertical and
%a horizontal symmetry axis, because then it would allow
% be symmetric with respect to
% a $180^\circ$ rotation, and this is impossible:
If  $n\ge4$, then
lines 1 and $n$ partition the plane into four sectors.
The crossing of lines 2 and $n-1$ is
in one of these four sectors, and the rotation moves it to
the opposite sector.

 \section{Computational results}
 \label{some-results}
As mentioned, going through all 12-point AOTs takes around 200 CPU
hours.
We also ran the program for 13 points on a
parallel compute-cluster~\cite{CURTA},
which took about 3200 CPU days of computing time.
The number~$h$ of hull vertices and the symmetry is already
computed
as part of the lex-min test; thus we might as well record these data.

For the purpose of illustration,
we decided to take some more statistics: the number of halving-lines
and the crossing number.

These data can be computed from the wiring-diagram:
The number of crossings at level $k$ in the wiring-diagram is the
number of lines through pairs of points that have $k$ points
\emph{below}
them, and hence it is clear that these numbers are related to the
$k$-edges and $k$-sets.
In particular, by counting the crossings at the different levels, one
immediately obtains the number of halving lines. By a % simple and
remarkably
 simple %and
formula
of
Lov\'{a}sz, Vesztergombi, Wagner, and Welzl
\cite{quad-k-sets-2004}, the number of crossings can be calculated
directly from the number of $k$-edges for all $k$.

Since we did not know what interesting phenomena might emerge from
the data, we decided
not to do any aggregation during the enumeration.
We maintain the number of AOTs for each combination of the characteristics ($n$, $h$, symmetry,
halving-lines, crossings), and in the end, we write the nonzero
numbers to a log-file, thus relieving the enumeration of the task to
make a statistical analysis. The result
file
with these raw data
is available in the
repository.\footnote{\cite{github}, file
  \texttt{results/crossing+halving-results-13.txt}; the data for the
  (realizable) OTs in Table~\ref{tab:mirror} are in the file
  \texttt{results/realizable11-10.txt}.}
All statistics reported in Sections~\ref{sec:hull}--\ref{sec:symmetries} were
extracted from these files by
 \textsc{Python} scripts.
\subsection{Number of convex hull points}
\label{sec:hull}
Table~\ref{tab:hull} counts the AOTs by size of the convex hull $h$.
This can be compared to Table~2 in~\cite{AICHHOLZER20072}, where the
corresponding data are given for (realizable) order types up to $n=11$.
\begin{table}[htb]
  \small
  \def\.{\hbox{\hskip-0,5pt '\hskip-0,6pt }}
  \def\.{\hbox{\hskip-0,3pt ,\hskip-0,3pt }}
  \centering
  \begin{tabular}{|@{ }l@{ }|@{ }r@{ \;}r@{ \;}r@{ \,}r@{ \;}r@{ \;}r@{ \;}r@{ }|}
    \hline
%    &$n=6$
    &$n=7$
    &$n=8$
    &$n=9$
    &$n=10$
    &$n=11$
    &$n=12$
    &$n=13$
    \\
    \hline
%   $h=3$ & 6& 49& 1\,178& 55\,239& 4\,879\,546& 786\,103\,220& 229\,258\,881\,954\\
% $h=4$ & 6& 59& 1\,468& 70\,482& 6\,324\,559& 1\,031\,019\,051& 303\,315\,298\,426\\
% $h=5$ & 3& 22& 570& 28\,234& 2\,630\,639& 440\,348\,013& 132\,120\,240\,798\\
% $h=6$ & 1& 4& 90& 4\,552& 450\,300& 79\,039\,502& 24\,562\,198\,935\\
% $h=7$ & & 1& 8& 311& 33\,969& 6\,447\,723& 2\,124\,883\,478\\
% $h=8$ & & & 1& 11& 1\,146& 241\,522& 87\,484\,087\\
% $h=9$ & & & & 1& 22& 4\,006& 1\,683\,531\\
% $h=10$ & & & & & 1& 33& 14\,410\\
% $h=11$ & & & & & & 1& 62\\
% $h=12$ & & & & & & & 1\\\hline
    %
$h=3$ & 49& 1\.178& 55\.239& 4\.879\.546& 786\.103\.220& 229\.258\.881\.954& 120\.410\.822\.315\.097\\
$h=4$ & 59& 1\.468& 70\.482& 6\.324\.559& 1\.031\.019\.051& 303\.315\.298\.426& 160\.356\.153\.417\.352\\
$h=5$ & 22& 570& 28\.234& 2\.630\.639& 440\.348\.013& 132\.120\.240\.798& 70\.900\.318\.730\.166\\
$h=6$ & 4& 90& 4\.552& 450\.300& 79\.039\.502& 24\.562\.198\.935& 13\.533\.084\.234\.118\\
$h=7$ & 1& 8& 311& 33\.969& 6\.447\.723& 2\.124\.883\.478& 1\.222\.365\.995\.348\\
$h=8$ & & 1& 11& 1\.146& 241\.522& 87\.484\.087& 53\.890\.715\.843\\
$h=9$ & & & 1& 22& 4\.006& 1\.683\.531& 1\.154\.715\.041\\
$h=10$ & & & & 1& 33& 14\.410& 11\.618\.261\\
$h=11$ & & & & & 1& 62& 51\.210\\
$h=12$ & & & & & & 1& 101\\
$h=13$ & & & & & & & 1\\\hline
    sum % $\Sigma$
    &135&3\.315&158\.830&14\.320\.182&2\.343\.203\.071&691\.470\.685\.682&366\.477\.801\.792\.538\\
    average $h$&
3.8815& 3.8793& 3.8935&
                        3.913\,29&3.928\,582~~&3.940\,299\,5\ ~&3.949\,367\,11\ ~\\
 % 3.88148148&3.87933635&3.8934773&3.91329021&3.92858207&3.94029949&3.94936711\\
    \eqref{eq:welz}&
  3.8182& 3.8621& 3.8919&
                          3.913\,04&3.928\,571~~&3.940\,298\,5\ ~&3.949\,367\,09\ ~\\
                     %3.81818182&3.86206897&3.89189189&3.91304348&3.92857143&3.94029851&3.94936709\\
%
%    sum %    $\Sigma$
%        &16&135&3\,315&158\,830&14\,320\,182&2\,343\,203\,071&691\,470\,685\,682\\
%    average $h$&
%    3.9375& 3.8815& 3.8793& 3.8935& 3.91329&3.928582&3.9402995\\
%\eqref{eq:welz}&   3.7500& 3.8182& 3.8621& 3.8919& 3.91304&3.928571&3.9402985\\
    \hline
  \end{tabular}
  \medskip
  \caption{Number of abstract order types with $n$ points in total and
    $h$ points on the convex hull.
The last row is the value of formula \eqref{eq:welz}, which has been
discussed in \autoref{sec:duplicates} on p.~\pageref{eq:welz}.
  }
  \label{tab:hull}
\end{table}

\begin{table}[htb]
  \centering
  \begin{tabular}{|l|rrrrrr|}
    \hline
& $h=3$ & $h=4$ & $h=5$ & $h=6$ & $h=7$ & $h=8$ \\\hline
$n=10$ & 0.340746 &  0.441654 &  0.183702 &  0.031445 &  0.002372 &  0.000080 \\
$n=11$ & 0.335482 &  0.440004 &  0.187926 &  0.033731 &  0.002752 &  0.000103 \\
$n=12$ & 0.331553 &  0.438652 &  0.191071 &  0.035522 &  0.003073 &  0.000127 \\
$n=13$ & 0.328562 &  0.437560 &  0.193464 &  0.036927 &  0.003335 &
                                                                    0.000147
    \\
\hline    
  \end{tabular}
  \caption{Relative frequencies of convex hull sizes $h$}
  \label{tab:relative-h}
% python3 makereport.py ../crossing+halving-results-13.txt
\end{table}

\begin{table}[htb]
  \centering
  \begin{tabular}{|rr|c|rr|c|rr|rr|rr|}
    \hline
    $X$&\#AOT&&
    $X$&\#AOT&&
    $X$&\llap{\#}AOT&
    $X$&\llap{\#}AOT&
               $X$&\llap{\#}AOT\\
    \hline
153&1&&250&9{,}599{,}727{,}792&&451&41&459&11&470&11\\
154&15&&251&9{,}774{,}280{,}765&&452&76&461&41&471&1\\
155&215&&252&10{,}813{,}519{,}833&&453&68&462&12&472&1\\
156&1{,}354&&253&10{,}549{,}648{,}258&&454&119&463&21&474&1\\
157&4{,}066&\smash{$\vdots$}&254&9{,}551{,}226{,}473&\smash{$\vdots$}&455&33&464&1&477&5\\
158&6{,}966&&255&9{,}720{,}622{,}387&&456&46&465&10&479&1\\
159&13{,}904&&256&10{,}543{,}935{,}293&&457&1&467&2&486&1\\
160&42{,}950&&257&10{,}332{,}151{,}661&&458&38&468&7&495&1\\
    \hline
\end{tabular}
  \smallskip
  \caption{The number of AOTs of 12 points with $X$ crossings, for
    selected values of $X$}
  \label{tab:crossings}
\end{table}

\begin{figure}[htbp]
  \centering
\hbox to \linewidth{\hss %\null  \hskip -8mm
  \includegraphics[width=1.05\textwidth
  ]{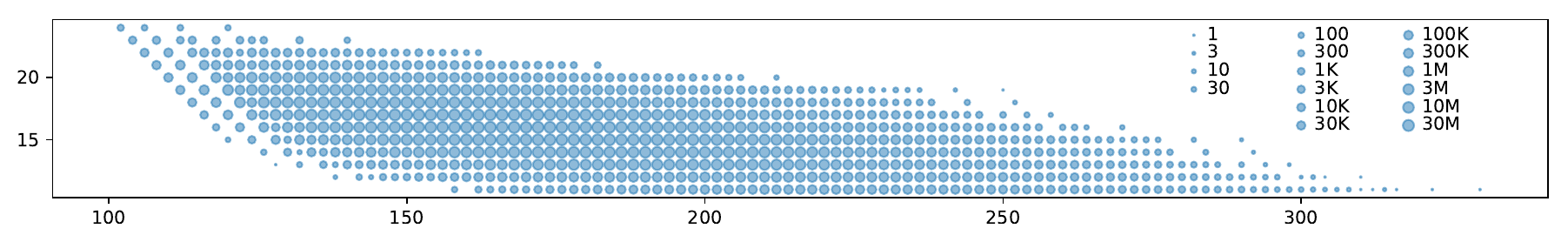}\hskip -2mm }
  \caption{Scatter-plot of crossing number (horizontal axis)
    versus number of halving-lines (vertical axis) for AOTs with
    $n=11$ points.
    The area of each dot represents the frequency, on a logarithmic
    scale. Thus the tails are really not as heavy as they appear.
    One can see that the crossing number and
    the number of
    halving-lines are negatively correlated.
    The crossing number ranges between 102 and 330, and it is always an
    even number.
    The number of
    halving-lines ranges between 11 and~24.
  }
  \label{fig:scatter11}
\end{figure}

Table~\ref{tab:relative-h} shows the relative frequencies of the
various convex hull sizes $h$, for the larger values $n= 10,11,12,13$. They
seem to converge to some limiting distribution. We are not aware of
any theoretical results that would predict the limiting
frequency of, say, triangular
convex hulls. This should be related %(but not by a direct proportion)
to the expected number of triangular faces in a ``random'' PSLA.

\subsection{Crossing numbers and halving-lines}
Table~\ref{tab:crossings} deals with the number of crossings in a
straight-line drawing of the complete graph.
We report results only for 12 points.
The smallest number of crossings is 153 (which has been known for a
long time), and is achieved by a unique AOT.
The largest number of crossings is $495=\binom{12}4$, and is achieved
by a unique AOT, namely by points in convex position.
The next-largest number of crossings is 486, and it is again
 achieved
 by a unique AOT.
 There are a few more gaps, as visible in the table.
 For every number $X$ in the range 153--459, there is at least one AOT with that
 number of crossings. The most frequent number of crossings is
 $X=252$; we can see that the frequencies do not vary monotonically but fluctuate up and down
 in the vicinity of this value.
%
%gaps:
%460, 466, 469, 473, 475, 476, 478, 480, 481, 482, 483, 484, 485, 487, 488, 489, 490, 491, 492, 493, 494

% halving-edges
% 4 [(2, 1), (3, 1)]
% 5 [(5, 1), (6, 1), (7, 1)]
% 6 [(3, 5), (4, 5), (5, 5), (6, 1)]
% 7 [(7, 5), (8, 13), (9, 43), (10, 44), (11, 27), (12, 3)]
% 8 [(4, 284), (5, 1018), (6, 1299), (7, 633), (8, 77), (9, 4)]
% 9 [(9, 673), (10, 4814), (11, 17215), (12, 36308), (13, 47038), (14, 35539), (15, 14423), (16, 2590), (17, 220), (18, 10)]
%10 [(5, 517423), (6, 2584235), (7, 4865400), (8, 4290426), (9, 1757011), (10, 283580), (11, 21389), (12, 713), (13, 5)]
%11 [(11, 2319053), (12, 20572280), (13, 88684444), (14, 241315734), (15, 449811080), (16, 582445724), (17, 515463082), (18, 301942617), (19, 111806473), (20, 24983674), (21, 3506778), (22, 330983), (23, 20435), (24, 714)]
%12 [(6, 11413112354), (7, 69595028879), (8, 171012193771), (9, 217728989695), (10, 152467416291), (11, 57295395849), (12, 10779353518), (13, 1112775181), (14, 64456370), (15, 1939706), (16, 23972), (17, 95), (18, 1)]

 Figure~\ref{fig:scatter11} shows the joint distribution of both parameters,
the   crossing number % (horizontal axis)
and the number of halving-lines. % (vertical axis) for AOTs with

\begin{table}[htb]
  \centering
  \vbox{
\halign{\strut\hfil$#$ &\ \hfil$#$ &\ \hfil$#$ &\ \hfil$#$  &\ \hfil$#$ 
  &\ \hfil$#$\cr
&\hbox{[\href{https://oeis.org/A006247}{A006247}]}
& \hbox{unsymmetric}
&\hbox{mirror-sym.}
&\hbox{rot.\,sym.}
&\hbox{[\href{https://oeis.org/A006246}{A006246}]}
 \cr 
 n & \hbox{(unoriented) AOTs} & \hbox{AOTs}
&  \hbox{AOTs}
&  \hbox{AOTs}
&  \hbox{oriented AOTs}
\cr
\noalign{\hrule}
3& 1 & 0 & 1 & 0 & 1 \cr
4& 2 & 0 & 2 & 0 & 2 \cr
5& 3 & 0 & 3 & 0 & 3 \cr
6& 16 & 4 & 12 & 0 & 20 \cr
7& 135 & 105 & 28 & 2 & 242 \cr
8& 3{,}315 & 3{,}085 & 225 & 5 & 6{,}405 \cr
9& 158{,}830 & 157{,}981 & 825 & 24 & 316{,}835 \cr
10& 14{,}320{,}182 & 14{,}306{,}748 & 13{,}103 & 331 & 28{,}627{,}261 \cr
11& 2{,}343{,}203{,}071 & 2{,}343{,}126{,}871 & 76{,}188 & 12 & 4{,}686{,}329{,}954 \cr
12& 691{,}470{,}685{,}682 & 691{,}468{,}293{,}616 & 2{,}358{,}635 & 33{,}431 & 1{,}382{,}939{,}012{,}729 \cr
13& 366{,}477{,}801{,}792{,}538 & 366{,}477{,}779{,}812{,}782 & 21{,}954{,}947 & 24{,}809 & 732{,}955{,}581{,}630{,}129 \cr
}}
\smallskip
\caption{AOTs with various symmetries.
 The column headings link to the
    corresponding entries of the
Online Encyclopedia of Integer
Sequences \cite{OEIS}.
 }
  \label{tab:symmetries}
\end{table}

\begin{table}[htb]
  \centering
  \begin{tabular}{|r|rrr@{\,}r@{\,}r@{\,}r@{\,}r@{\,}r@{\,}r@{\,}r@{\,}r@{\,}r@{\,}r|r@{
    }r@{\
    }r@{\ }r@{\ }r|}\hline
$n$&$D_{1}$ & $D_{2}$ & $D_{3}$ & $D_{4}$ & $D_{5}$ & $D_{6}$ &
                                                                $D_{7}$
    & $D_{8}$ & $D_{9}$ & $D_{10}$ & $D_{11}$ & $D_{12}$ & $D_{13}$&
   $C_2$ & $C_3$ & $C_4$ & $C_5$ & $C_6$    \\
\hline    
3&  &  & 1 &  &  &  &  &  &  &  &  &  &  & &&&&\\
4&  &  & 1& 1&  &  &  &  &  &  &  &  &  &  & & & &  \\
5& 2&  &  &  & 1&  &  &  &  &  &  &  &  &  & & & &  \\
6& 7& 1& 2&  & 1& 1&  &  &  &  &  &  &  &  & & & &  \\
7& 26&  & 1&  &  &  & 1&  &  &  &  &  &  &  &2& & &  \\
8& 218& 4&  & 1&  &  & 1& 1&  &  &  &  &  & 4& &1& &  \\
9& 818&  & 6&  &  &  &  &  & 1&  &  &  &  &  &24& & &  \\
10& 13{,}059& 27& 11&  & 4&  &  &  & 1& 1&  &  &  & 234&93& &4&  \\
11& 76{,}186&  &  &  & 1&  &  &  &  &  & 1&  &  &  & & &12&  \\
12& 2{,}358{,}210& 303& 111& 7&  & 2&  &  &  &  & 1& 1&  & 29{,}573&3{,}765&86& &7 \\
13& 21{,}954{,}912&  & 34&  &  &  &  &  &  &  &  &  & 1&  &24{,}809& & &  \\
    \hline
    \multicolumn{19}{|c|}{(realizable) OTs}\\
    \hline
9& 818&  & 6&  &  &  &  &  & 1&  &  &&&  &23\rlap{$^*$}& &  &\\
10& 13{,}021\rlap{$^*$}& 27& 10\rlap{$^*$}&  & 4&  &  &  & 1& 1&&&  &
        234&84\rlap{$^*$}& &3\rlap{$^*$}
   & \\
11& 75{,}802\rlap{$^*$}&  &  &  & 1&  &  &  &  &  & 1&&&  & & &7\rlap{$^*$}& \\
    \hline
  \end{tabular}
  \medskip
  \caption{The symmetric AOTs according to their symmetry group. The
    last three rows concern OTs, for those cases where the set of
     OTs is known ($n\le 11$) and differs from the set of AOTs
     ($n\ge9$).
The differences to AOTs are marked. The majority of
the difference set (column~$\Delta$ in Table~\ref{tab:number})
belongs to the class $C_1$ with no symmetries at all, which is
not shown in this table.
   }
  \label{tab:mirror}
\end{table}

\subsection{Symmetries}
\label{sec:symmetries}

As mentioned in \autoref{sec:duplicates}, our algorithm delivers the symmetries of
an AOT for free, as part of the lex-min test that is necessary to pick
a single PSLA among the several PSLAs representing the AOT.
Table~\ref{tab:symmetries} classifies the AOTs (first column) according to the types
of symmetries that they have. The second column gives the AOTs that
have no symmetry at all, and these are the vast majority.
The third column counts AOTs that have a mirror symmetry
(possibly including a rotational symmetry as well).
The fourth column gives the AOTs that have a non-trivial symmetry that
is purely rotational
 (without mirror symmetry).
The first column is the sum of columns 2--4.

A mirror symmetry will \emph{reverse} all orientations, and thus,
there
can be different opinions whether it should be regarded as a symmetry operation.
The last column is the number of \emph{oriented AOTs}, where an AOT
and its mirror are counted as distinct objects (unless the AOT is
mirror-symmetric).
It is obtained by taking columns 2 and 4 twice and adding column~3 once.

Table~\ref{tab:mirror} gives a more refined account of columns 3 and
4 of 
Table~\ref{tab:symmetries}, classifying AOTs by their symmetry groups.
We use the same notations $C_k$ and $D_k$ as for the symmetry groups of
finite objects in the plane (the cyclic and dihedral groups), although
our groups act in a purely combinatorial way on AOTs, by permuting
the points. Since such a symmetry must preserve the convex hull
vertices and the adjacency between them,
it %the symmetry group of an AOT
must be isomorphic to a subgroup of the symmetry group of a regular $h$-gon, if there are
$h$ hull vertices.

$D_k$ is the symmetry group of a regular
$k$-gon: the \emph{dihedral group} of order $2k$. For each~$n$, we have one
AOT with symmetry group $D_n$, namely convex position, which
corresponds to the regular $n$-gon in the geometric setting. In
addition, if $n$ is even, we can have an $(n-1)$-gon with a point in
the center, having symmetry group $D_{n-1}$. The most frequent group
is the group $D_1$, which has a single mirror-symmetry as the only
nontrivial element.

$C_k$ is the rotational symmetry group of a regular $k$-gon, i.e., a $k$-fold rotation.
$C_2$~corresponds to a rotation by $180^\circ$, or equivalently, a
reflection in a central point.

We see that many fields in the table are empty.
There are systematic reasons for this.
For example, if a set of
$n$ points has a rotational symmetry of order 3 ($C_3$, or any of its
supergroups $D_3$, $C_6$, $D_6$, $D_9$, or $D_{12}$), then $n$ must be a
multiple of 3 or a multiple of 3 plus~1 (with a fixpoint in the
``center''),
cf.~\cite[Theorem 1.5]{goaoc-welzl-2023}.
A $C_2$ symmetry cannot exist for odd~$n$, because it
would
have a fixpoint
in the center, and ``opposite'' points would have to be aligned
with the center, which is excluded in a simple AOT.
% \footnote
% {Goaoc and Welzl \cite[Theorem 1.5]{goaoc-welzl-2023}
% claim that the symmetry group of an order type must be isomorphic to a
% cyclic group, thus mistakenly
% excluding mirror symmetries.
% ($D_1$ is indeed isomorphic to a cyclic group, but the other dihedral groups
% are not.)
% Moreover, the
% example that they show for a symmetry group of order 2
%  \cite[Figure~4 left]{goaoc-welzl-2023} actually has the symmetry group~$D_2$.
% }

\paragraph{Counting of AOTs by enumerating symmetric AOTs.}
There is a relation between the number of AOTs
with prescribed symmetries and the number of PSLAs:
Each AOT corresponds to a certain number of PSLAs, depending
on the symmetry group. Thus if we know the entries in
Table~\ref{tab:mirror}, together with the unsymmetric AOTs in
the second column of Table~\ref{tab:symmetries}, we can work out the
number of PSLAs.
(This is actually how the program computes the correct number of PSLAs
even though it prunes branches of the enumeration tree and does not
visit each PSLA individually.)

We could use this relation in the other direction.
We might think about counting the various symmetric AOTs for $n>13$ by
enumerating them directly, since their number should still be
manageable. Together with the number of PSLAs, which is known up to
16 lines, we can then calculate the number of non-symmetric AOTs, and hence the
total number of AOTs.

\subsection{The number of cutpaths of a PSLA}

\begin{figure}[htb]
  \centering
%  \noindent  \hskip -5mm
%  \includegraphics[width=1.08\textwidth]{plot-cutpaths-8}
  
    \null\hskip -17,6mm
  \includegraphics[width=1.103\textwidth]{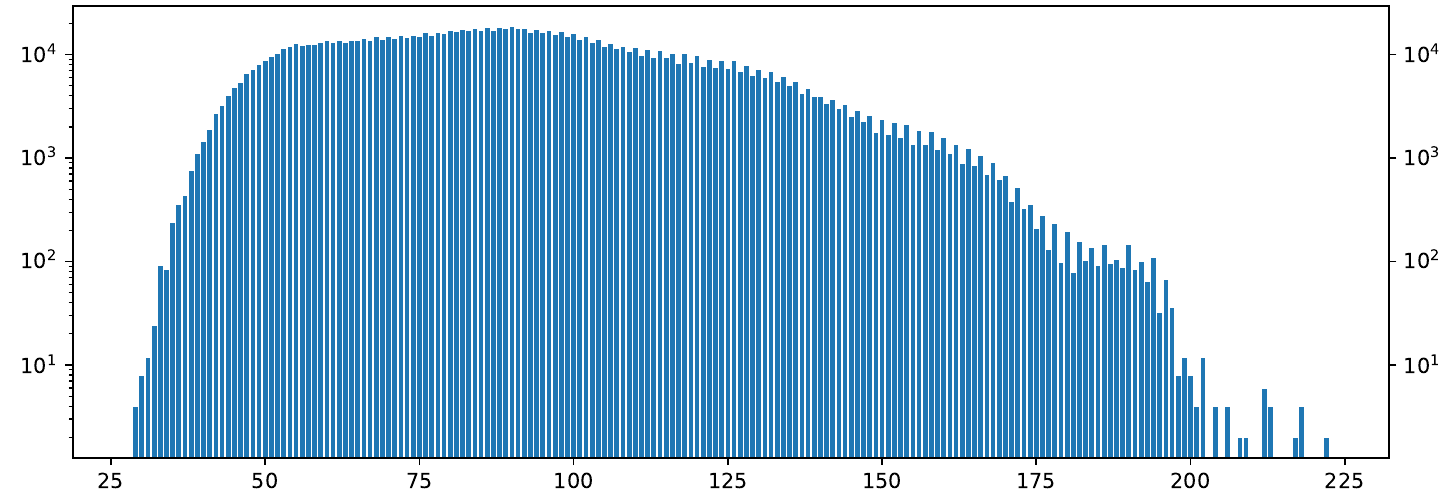}%
  \hskip -17mm\null
  \caption{The distribution of the number of cutpaths of the
1,232,944
    PSLAs with
    8~pseudolines.
    The frequencies are on a logarithmic scale.
    For symmetry reasons,
    every number of cutpaths occurs with an even frequency
    (see \autoref{sec:symmetry-PSLA}).
    % %Thus, the reported frequencies have been divided by 2.
    The number of cutpaths ranges between 29 and 222, and the average
    is 90.85, which equals
 the quotient of the number of PSLAs with 9 and with 8
    pseudolines,
    see Table~\ref{tab:number}.
% 112018190/1232944 ~ 90.85423993303831  
  }
  \label{fig:cutpaths}
\end{figure}

An extra pseudoline that runs from the bottom face to the top face in  a given PSLA is
called a {cutpath} (see Figure~\ref{fig:dualDAG}).
The number of cutpaths is equal to the number of children of the
corresponding node in the enumeration tree.
\autoref{fig:cutpaths} shows the distribution of the number of
cutpaths
for the PSLAs with 8 pseudolines.\footnote{
Frequencies of the number of children have been computed for PSLAs
with up to 10 pseudolines;
  the data
are available in the
repository
\cite[%file
  \texttt{results/cutpaths+grandchildren-results-10.txt}]{github}.}%
\footnote
{We remark that Knuth~\cite{axioms-hulls}
has computed some statistics about the number of cutpaths,
in particular about the minimum, the 
maximum and the average~\cite[Table at the bottom of p.~39]{axioms-hulls}.
However, he counts a different class of cutpaths, namely
the paths from the top face $T$ to the bottom face $B$ in
 the spherical arrangement, see
 \autoref{fig:spherical-model} and
 \cite[Example (9.9) on p.~38]{axioms-hulls}.
%consists of two justaposed copies of the same \psla.
%Such an arrangement arises if the left half and the right half of
%five roundabout connections from the left edge to the right edges.
%The second copy is vertically reflected; each pair of curves intersects twice, and
%the order of the \psl s at the right edge is the same as the order at the left edge.
(Note that the %Knuth's
conjecture
%that Knuth expressed
about the maximum number $n2^{n-2}$ of cutpaths
that was expressed in this context
% \cite[ ... Example (9.9) on p.~38]{axioms-hulls},
has
been disproved by
Ond\v rej B\'\i lka; see 
\cite{felsner-valtr-11,d-ibnpa-25} for details.)
}

  \section{Estimating the number of PSLAs with larger \texorpdfstring{$n$}{n}}
\label{sec:estimate}

  Knuth~\cite{knuth75_estim}
 observed that there is an easy method to obtain
  an unbiased estimate for the
  number of leaves of a tree, by
  %performing
  making
  a random ``dive'' into the depth of the tree:
%can be obtained by
%following the path
Iteratively proceed to a random child and
  multiply the encountered vertex degrees,
see also
%TAOCP
%Volume 4B
%Combinatorial Algorithms, Part 2
\cite[Sect.~7.2.2, pp.~46--51,
Corollary~E]
{kn5}.
It is straightforward to prove that the expected value of this
estimate is indeed the
number of leaves of the tree.

In this way we can estimate the size of the deeper levels of the
enumeration tree, beyond the range that can be explicitly enumerated
one by one.
The procedure can be easily adapted to estimate the number of nodes
on \emph{each level} of the tree, up to a maximum depth,
or to give an unbiased estimate for the number of OAOTs and
AOTs.
%For the latter estimation, one simply has to 

For each node we need to compute
the number
of cutpaths (= the number of children),
and, except at the lowest level, pick a random child.
The number $c$ of cutpaths
equals the number of source-to-sink paths in the
dual DAG of the PSLA (see the
right part of \autoref{fig:dualDAG}).
It can be computed by a traversal of the graph in topological order
in linear time in the size of the graph, i.e., in $O(n^2)$ time.
In the same time, we can also, for a random number $r$ between 1 and $c$,
determine the $r$-th cutpath in lexicographic order. Thus, we can
determine a random cutpath in $O(n^2)$ time.
%(The saving of this method is not very spectacular. For
%$n=8$ pseudolines, the dual DAG has $1+\binom{n+1}2=37$ nodes
%{fig:cutpaths}

\begin{figure}[htb]
  \centering
  \noindent
  \null \hskip-5mm plus 1mm \null 
  \includegraphics[height=41.5mm]{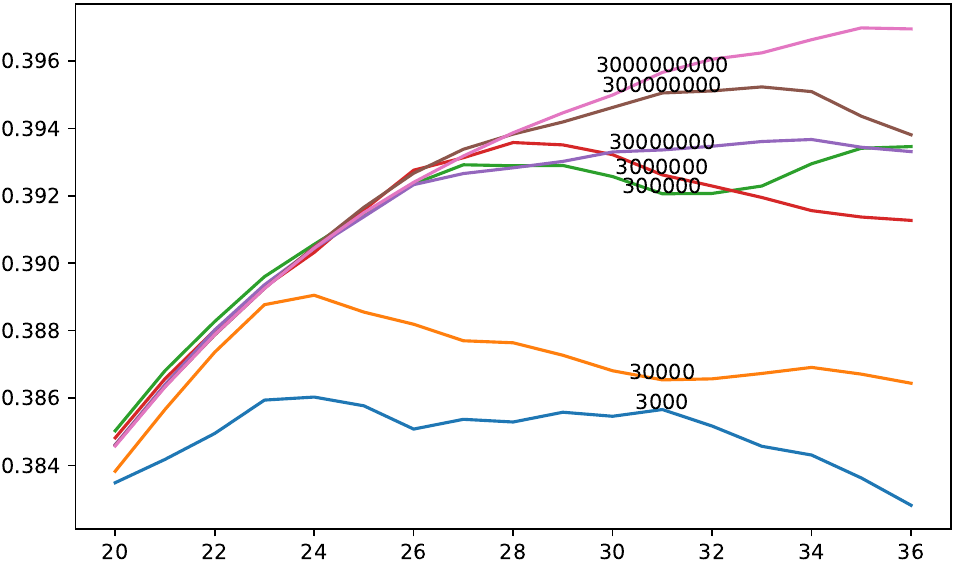}\hfill
  \includegraphics[height=41.5mm]{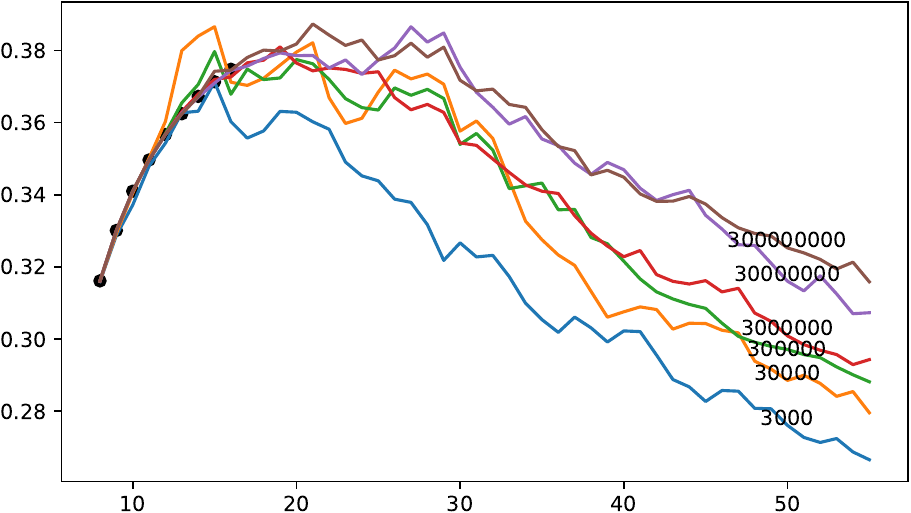}%
%  \null\hskip-0mm  \null
  \caption{%Estimating % the constant
    $(\log_2 \bar B
%    \textrm{\#PSLA}
    _n)/n^2$,
    where $\bar B_n$ is the estimate for \textrm{\#PSLA}$_n$ obtained as the
    average of $m$~samples, for various values of $m$ between 3000 and $3$~billion. %\times 10^{9}$.
    The horizontal axis is the number $n$ of \psl s. The right graph
    uses a simpler and faster method.
  }
  \label{fig:random}
\end{figure}

\autoref{fig:random} shows the results of two experiments
to estimate \#PSLA in this way.
It is known that the number of pseudoline arrangements grows
asymptotically like
$2^{\Theta(n^2)}$;
%const$^{n^2}$;
more precisely,
\begin{displaymath}
   0.2721
  \le
  \liminf_{n\to\infty}
\frac{  
\log_2 %B
    \textrm{\#PSLA}
    (n)}{n^2}
  \le
    \limsup_{n\to\infty}
\frac{  
\log_2 %B
    \textrm{\#PSLA}
    (n)}{n^2}
\le  0.6496, %0.6571,
\end{displaymath}
and it is believed that  $(\log_2 %B
    \textrm{\#PSLA}
    (n))/n^2$ converges to a constant.
%    Both the
The best known upper bound (due to Dallant~% is due to Felsner and Valtr from 2011
\cite%{felsner-valtr-11}.
{d-ibnpa-25})
and lower bound~\cite{corteskuhnast.SoCG.2024}
 are quite recent and are still being improved.
%Therefore, in the plots, we have taken the logarithm of
%the estimate
%$\textrm{\#PSLA}_n$ and divided by $n^2$.

The left part of
\autoref{fig:random}
shows the estimate of
the number of PSLAs with $n$ lines,
taking the average $\bar B_n$ of the results from $3\times 10^9$
random paths %dives
to depth 36, % into the enumeration tree,
together with some partial estimates
obtained along the way from smaller subsamples.
In accordance with the exponential growth of \textrm{\#PSLA}, we have
plotted
the quantities
$(\log_2 \bar B
%    \textrm{\#PSLA}
_n)/n^2$.

This method, as currently implemented,
reaches its limits somewhere near %beyond
$n=70$. When
the number of children exceeds $2^{64}$, calculating
the number of children requires multi-precision arithmetic,
and
a more elaborate random number generator would
be needed to select a random child.
We set the cutoff, beyond which no reasonable random choice
can be expected, at a maximum of $2^{55}$ children.
This limit was often exceeded for large samples already for $n=69$.

We therefore considered a
 simplified sampling procedure, whose results are shown in
the right part of
\autoref{fig:random}.
When incrementally inserting the $n$-th pseudoline, we always pick a
random edge on the upper boundary of the current face where the
pseudoline
should leave the face (see \autoref{fig:explore-face}).
This procedure generates a different enumeration tree for the PSLAs,
which is very much
refined in comparison to the tree of \autoref{fig:tree},
and it is simpler and much faster:
The $n$-th pseudoline can be inserted with $n-1$ consecutive random decisions
in $O(n)$ time.
(By the zone theorem, the total number of edges of the faces that a
cutpath visits is $O(n)$.)
This allowed us to go up to $n=200$ pseudolines,
taking
% another run with
300 million
samples. %experiments.
\autoref{fig:random} shows the estimates up to
$n=55$. The black dots represent the true values of
\textrm{\#PSLA}$(n)$, which are known up to $n=16$.
Comparing to the original experiment in the left graph, we see that
the fluctuations are much larger and occur earlier.
%Also, the estimate for the original experiment reaches values beyond
%0.395 for $m= 300$ million or larger, whereas the simplified
%experiment stays well below 0.39. NAJA.

 Files with the complete results, including also estimates of AOTs and OAOTs, are
available in the repository.\footnote{\cite[%{github}, files
  \texttt{results/random-estimate-results-36.txt},
  % and
  \texttt{results/random-face-estimate-results-200.txt}]{github}}

In both experiments,
we can notice %see %the tendency to
a systematic underestimation,
with a few exceptional overestimations.
%as $n$ gets larger.
The reason is that the estimate, while having
 the correct expectation,
 has an extremely large variation. Most of the subtrees at a given level are
 smaller than average,
and very few subtrees are huge.
With a uniform choice
of a child, it is unlikely that
 the rare huge subtrees are hit, and
 as long as none of these few % rare huge
 subtrees is hit,
the estimate
 is too low.
 The same behavior has been observed in other contexts, for example,
 when estimating the cost %size
 of backtracking trees~\cite{knuth75_estim} and branch-and-bound trees.

\begin{figure}[htbp]
  \centering
  \noindent  \hskip -5mm
%  pdfcrop --margins 10 ~/Figure_1.pdf plot-7-children-grandchildren.pdf
  \includegraphics[width=0.91\textwidth]{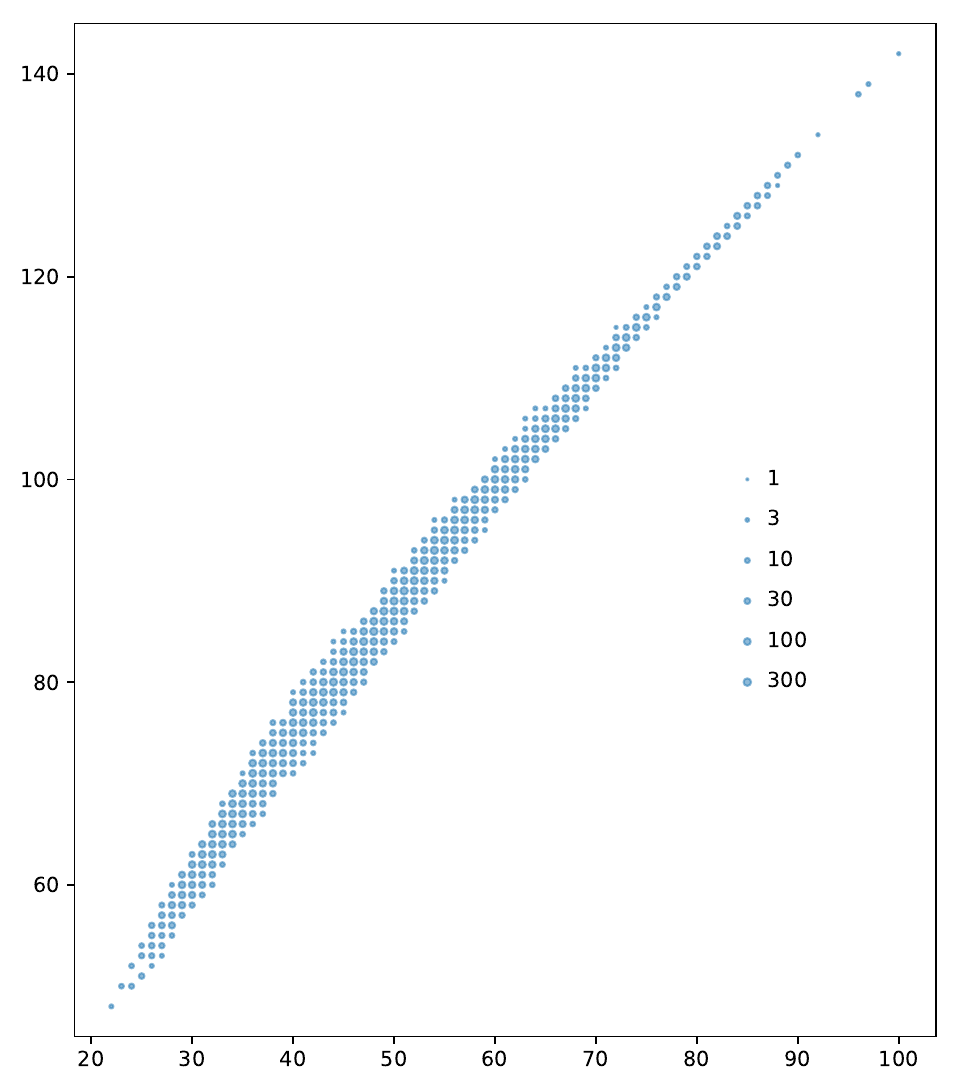}
  \caption{Scatter-plot of children versus grandchildren, for PSLAs with 7 pseudolines.
    The horizontal axis gives
the number of cutpaths of each PSLA, or in other words,
    the degree of each node (the number of
children) in the enumeration tree.
The vertical axis gives the average degree of those children, or in
other words, the number of grandchildren divided by the number of
children,
rounded to the nearest integer.
    The area of each dot represents the frequency, on a logarithmic
    scale.
The number of cutpaths ranges between 22 and 100.
  }
  \label{fig:grandchildren}
\end{figure}

 One can try to counteract this phenomenon % effect
 by deviating
 from the uniform selection among the children,
 favoring the children that have themselves many children.
 The estimation formula can accommodate 
 any nonuniform selection: Instead of the product of
 the degrees, one takes the inverse of the product of the
 probabilities
 of the chosen children.

 In fact, it seems that the children of high-degree  vertices tend to
 have
 % children with
 high degree themselves.
 \autoref{fig:grandchildren} plots cutpaths of 7-line PSLAs in
relation to the average number of
cutpaths of their 8-line children.
The correlation between the degrees of vertices and
the degrees of their children is clearly visible.
%One clearly sees that PSLAs with many cutpaths tend to have children
%with many cutpaths.
For PSLAs with 8 or 9 pseudolines, the diagram looks qualitatively the
same.\footnote{The raw data
are available in the
repository
\cite[%file
  \texttt{results/cutpaths+grandchildren-results-10.txt}]{github}.}
This phenomenon exacerbates the poor behavior of a uniform choice
of a child.
 One should deviate from this uniform selection among
the children in a stronger way than just by weighing them by
\emph{their} number of children.
Different strategies in this direction might be evaluated experimentally.

 Another experiment that might be worth trying would be
 to start
 a single random path at each of the 18 billion nodes at level 10,
 as opposed to starting, say, 18 billion random paths from the root.
 It would be interesting to see whether % how much
 this improves the estimates.
Other methods of improving the estimate for the tree size have been proposed, for
example \emph{heuristic sampling}, due to Pang C. Chen
\cite{Chen-sampling-1992}
(see also \cite[Section 7.2.2.9]{kn4c-draft}),
or \emph{stochastic enumeration}, due to
Reuven Rubinstein~%
% Rubinstein 2013
\cite{Rubinstein-2013},
% for ``counting'' trees 2017
see also~\cite{Vaisman-Kroese-2017-trees}.
There is a general framework
 that encompasses these variations~\cite{r-ennpd-26}.

\section{Extensions}

There are many ways in which one could think of extending the program.
\begin{itemize}
  
\item In this article, we have concentrated on AOTs. 
  We have used PSLAs mainly as a tool to enumerate AOTs, but
  PSLAs can also be considered in their own right. They might be
  counted or classified with respect to different criteria, like
  projective equivalence classes or affine equivalence classes
  (cf.\  Figure~\ref{fig:PSLA}).
  % \item
  Each of these classes could be further classified by their symmetries,
  as in
\autoref{sec:symmetries}.
\item It should not be hard to generate ``partial'' \psla s, in which lines are not forced to cross;
  see Figure~\ref{fig:partialPSLA} for an example.

\begin{figure}[htbp]
  \centering
  \includegraphics[scale=0.62]{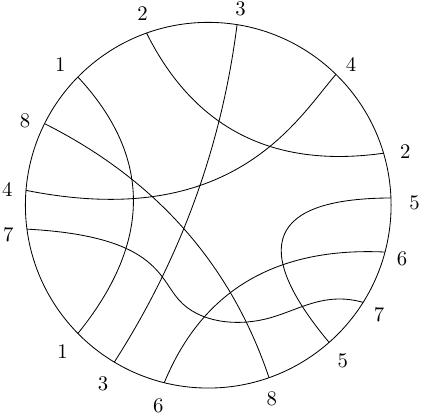}
  \caption{Example of a partial PSLA}
  \label{fig:partialPSLA}
\end{figure}

\item Random generation. It is easy to generate a random PSLA by
following a random path into the tree. This will, however, not be a uniform
selection.

To generate uniformly random \psla s with a particular number $n$ of \psl s, say $n=11$,
one could, in a preprocessing stage, compute and store the size of all 112
million subtrees at level~9. This allows to choose a random subtree
according to the correct distribution
quickly.
A~random leaf in that subtree then can be computed by exploring
the subtree, which contains on average about 50,000 leaves.
In going through these leaves, one does not need to go through every
leaf on the last
level (level~11); it mostly suffices to determine the number of children at the
previous level (level~10).
This method should therefore be reasonably fast.
A few
particularly large subtrees could be further preprocessed to speed up
the process.
The method can be adapted to generate a random AOT.

\item Non-simple \psla s, in which more than two \psl s are allowed to
  cross in a point, are challenging.
  In the language of oriented matroids, they correspond to nonuniform
  oriented matroids, and they have been enumerated by a method of
  Finschi and Fukuda~\cite{dcg-FinschiFukuda02}, also in higher dimensions,
see \cite{OM-homepage} for a catalog.
  These are much more numerous, see also
Table~1 in  \cite{dcg-FukudaMM13},
  % for an overview,
  where also the
  realizability % question
  is considered.
 Handling them by our approach would
  involve a
  redesign of the data structures from scratch.

  \iffalse
  A006247:  
1, 1, 1, 2, 3, 16, 135, 3315, 158830, 14320182, 2343203071,
691470685682, 366477801792538 = \#AOT = Knuth $D_n$

Also the number of nonisomorphic [``nondegenerate''=] uniform acyclic rank 3 oriented
matroids on $n$ elements.
$D_n$ = number of topologically distinct, simple arrangements of pseudolines with a marked
cell, as discussed by Goodman and Pollack [30].

A063800: 1, 2, 4, 17, 143, 4890, 461053, 95052532

Number of nonisomorphic oriented matroids with n points in 2
dimensions (rank 3), see
 \cite[Table~1]{dcg-FukudaMM13}.

A006248: 1, 1, 1, 1, 1, 4, 11, 135, 4382, 312356, 41848591,
10320613331

Number of projective pseudo order types: simple arrangements of
pseudo-lines in the projective plane. = Knuth $E_n$
= simple uniform oriented matroids (orientation classes) according to
  \cite[Table~1]{dcg-FukudaMM13}.

\# PSLA = Knuth $B_n$  

A006246 = 1, 1, 1, 2, 3, 20, 242, 6405, 316835, 28627261, 4686329954
= oriented AOTs = Knuth $C_n$ !! Extend OEIS

Number of simple arrangements of pseudolines in the projective plane
with an oriented marked cell; number of oriented abstract order types
of $n$ points (distinguishing mirror-symmetric copies).

subsection{Relation to oriented matroids} ....

relabeling classes, isomorphism classes

oriented matroids
of rank 3.

\fi

\item A side issue are nice drawings of \psla s. The wiring diagram is
  simple to obtain but it is very jagged.
  (Another representation of \psla s via rhombus tilings is even more jagged.)
  Stretchability of \psl s to straight lines, on the other hand, is a % very
  hard problem,
  despite the relative
  scarcity of non-stretchable PSLAs
  % (and  non-realizable AOTs)
  of moderate size.
  As the number %$n$
  of \psl s increases,
  non-stretchable arrangements become more and more dominating.
  Moreover, when drawing \psla s
  by  straight lines such that
  all vertices of the
  arrangement are visible, one often  cannot avoid having
  some vertices crammed close together in a tiny area. 
  An algorithm % program
  to construct a drawing in which the \psl s don't
bend too much
while keeping good distances between the features
 would be nice to have. %an interesting task.
%  (Maybe it
%  would be an idea for a  
%Geometric Optimization Challenge\footnote
%%Computational Geometry
%{\url{https://cgshop.ibr.cs.tu-bs.de/}}, perhaps in connection with
%the random generation method mentioned above.)

\end{itemize}

%\paragraph{Acknowledgements.}

\bibliography{../ENUM}

%\clearpage
\appendix

\section{Benchmark comparison with the order-type database}
\label{sec:benchmark}
We compared the usage of the order-type database against our
enumeration approach, and
we found that generation from scratch can actually compete in terms of runtime.
We compared the following two tasks:
\begin{enumerate}
\item
[A.]
Read the 14,309,547 order types of 10 points from the database
and compute the size of the convex hull.
The convex hull can be computed in linear time,
since the first point is always a convex hull vertex and the other
points are ordered clockwise around this point.
The coordinates are 16-bit unsigned integers, and the orientation test
is performed by determinant computation in the usual way, with two
multiplications, using 64-bit integer arithmetic.
\item
[B.]  Generate the
  14{,}320{,}182 abstract order types of 10 points
by \texttt{NumPSLA} and report
  the size of the convex hull.
  The size of the convex hull is computed anyway as part
  of the lexicographic normalization procedure; thus it
  does not cost any extra runtime.

  (With the \texttt{-exclude} option, we could also
  generate the
  14,309,547 \emph{realizable} abstract order types,  
  but this did not make a noticeable difference in runtime.)
  
\end{enumerate}
Both programs
took about 10--20 seconds
with a slight advantage for one or the other program,
depending on the computer.

For task A, typically about 60\,\% of the total time %10 seconds
were
``system time'', for reading the file, and
40\,\% were
``user time'', for the actual computation.

The usual goal is to perform some more time-consuming checks or
calculations on each order type. Thus, the time for either reading the point
set from the file or for generating it is a minor issue.

%user	0m3.897s
%sys	0m5.569s

%\clearpage
\section{A Python version of the basic enumeration algorithm}
\label{python}

\begin{figure}[htb]
  \centering
  \includegraphics[scale=0.85
  ]{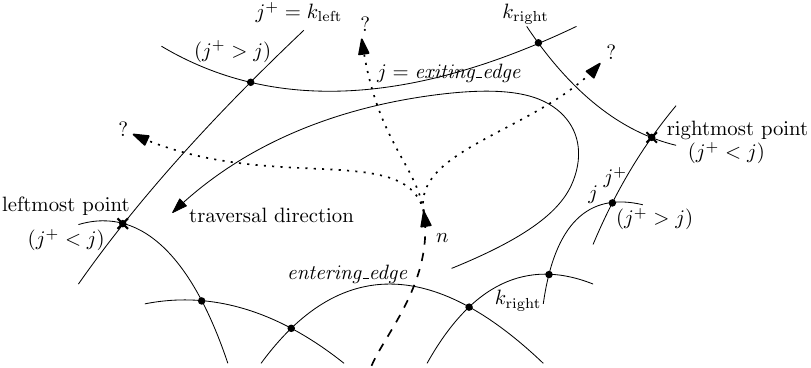}
  \caption{Threading line $n$ through a face}
  \label{fig:threading}
\end{figure}

The following program will carry out the basic enumeration of PSLAs.
The function \texttt{recursive\_generate\_PSLA\_start} is the outer
recursion, inserting the next pseudoline.
The function \verb|recursive_generate_PSLA| is the inner
recursion, extending pseudoline~$n$ into the next face by crossing an
edge.
Figure~\ref{fig:threading}
is a refined version of Figure~\ref{fig:explore-face},
illustrating the meaning of the variables
in the inner loop.

The program is available in the repository under the filename
\texttt{python/NumPSLA-basic.py}.
Starting the program with
\begin{quote}
  \texttt{python python/NumPSLA-basic.py 7}
\end{quote}
will count all $x$-monotone pseudoline arrangements with at most 7
lines
by running through each of them individually.
By importing the module \texttt{python/wiring\_diagram.py}, one can for example modify
the program to print wiring diagrams of all arrangements.

\ifminted
{\small % PLAY WITH THIS TO GET ON ONE PAGE
 % \large
%{%\captionsetup{type=listing}
%\begin{listing}
%  \inputminted{python}{../python/NumPSLA-basic.py}
\begin{minted}{python}
"The basic framework of NumPSLA, python version"
import sys
# from wiring_diagram import print_wiring_diagram #, IPE_end

def LINK(j, k1,k2): # make crossings with k1 and k2 adjacent on line j
    SUCC[j,k1] = k2
    PRED[j,k2] = k1
    
def Process_PSLA(n): # insert your code for processing the PSLA here:
    countPSLA[n] += 1
    # print(n,  countPSLA[n], ".".join(str(x) for x in localCountPSLA[3:n+1]))
    # print_wiring_diagram(n, SUCC, ipe=False)
                
def recursive_generate_PSLA(entering_edge, k_right, n):
    j,jplus = entering_edge,k_right
    while jplus>j: # find right vertex of the current cell F
        j,jplus = jplus,SUCC[jplus,j]
    # the right vertex is at the crossing of j and jplus
    if jplus==0: # F is unbounded
        if j==n-1: # F is the top face.
            LINK(n, entering_edge,0) # complete the insertion of line n
            localCountPSLA[n]+=1
            Process_PSLA(n)
            if n<n_max:
                localCountPSLA[n+1]=0 # reset child counter
                recursive_generate_PSLA_start(n+1) # thread the next pseudoline
            return;
        else: # jump to the upper bounding ray of F
            jplus=j+1; j = 0;
    while True:
        # scan the upper edges of F from right to left and try them out.
        k_right = j;
        j = exiting_edge = jplus;
        k_left = jplus = PRED[j,k_right];
        LINK(exiting_edge, k_left,n); # prepare for the recursive call
        LINK(exiting_edge, n,k_right);
        LINK(n, entering_edge,exiting_edge);
        
        recursive_generate_PSLA(exiting_edge, k_right, n) # enter the recursion
        
        LINK(exiting_edge, k_left,k_right); # undo the changes
        if jplus <= j: return
        #terminate at left endpoint of the face F or at unbounded ray (jplus=0)

def recursive_generate_PSLA_start(n):
    LINK(0, n-1,n);
    LINK(0, n,1); # insert crossing with line n on line 0
    recursive_generate_PSLA(0, 0, n);
    LINK(0, n-1,1); # undo the insertion of line n
          
n_max = int(sys.argv[1])

# Start the generation proper: 
PRED = {}; SUCC = {} 
LINK(1, 0,0);
LINK(0, 1,1);
  
countPSLA = [0]*(n_max+1)
localCountPSLA = [0]*(n_max+1)
recursive_generate_PSLA_start(2)
# IPE_end() # finish and close ipe-file, in case it was used.
print ("Number of PSLAs:", *countPSLA[2:])
\end{minted}
%  \caption{xxx}
  % \end{listing}
%  }
}
\fi

\end{document}